\documentclass[12pt]{article}
\usepackage{graphicx, amsmath, mathrsfs}
 \usepackage{color}
\usepackage{amsfonts}
\usepackage{verbatim}
\usepackage{latexsym}
\usepackage{enumitem}
\usepackage[toc,page]{appendix}

\makeatletter
\@addtoreset{equation}{section}

\begin{document}

\newcommand{\E}{\mathbb{E}}
\newcommand{\PP}{\mathbb{P}}
\newcommand{\RR}{\mathbb R}
\newcommand{\LL}{\mathbb{L}}

\newtheorem{theorem}{Theorem}[section]
\newtheorem{remark}{Remark}[section]
\newtheorem{lemma}{Lemma}[section]
\newtheorem{coro}{Corollary}[section]
\newtheorem{defn}{Definition}[section]
\newtheorem{assp}{Assumption}[section]
\newtheorem{expl}{Example}[section]
\newtheorem{prop}{Proposition}[section]

\newcommand\tq{{\scriptstyle{3\over 4 }\scriptstyle}}
\newcommand\qua{{\scriptstyle{1\over 4 }\scriptstyle}}
\newcommand\hf{{\textstyle{1\over 2 }\displaystyle}}
\newcommand\hhf{{\scriptstyle{1\over 2 }\scriptstyle}}

\newcommand{\eproof}{\indent\vrule height6pt width4pt depth1pt\hfil\par\medbreak}

\def\a{\alpha} \def\g{\gamma}
\def\e{\varepsilon} \def\z{\zeta} \def\y{\eta} \def\o{\theta}
\def\vo{\vartheta} \def\k{\kappa} \def\l{\lambda} \def\m{\mu} \def\n{\nu}
\def\x{\xi}  \def\r{\rho} \def\s{\sigma}
\def\p{\phi} \def\f{\varphi}   \def\w{\omega}
\def\q{\surd} \def\i{\bot} \def\h{\forall} \def\j{\emptyset}

\def\be{\beta} \def\de{\delta} \def\up{\upsilon} \def\eq{\equiv}
\def\ve{\vee} \def\we{\wedge}
\def\vv{\varepsilon}

\def\F{{\cal F}} \def\P{{\cal P}}
\def\T{\tau} \def\G{\Gamma}  \def\D{\Delta} \def\O{\Theta} \def\L{\Lambda}
\def\X{\Xi} \def\S{\Sigma} \def\W{\Omega}
\def\M{\partial} \def\N{\nabla} \def\Ex{\exists} \def\K{\times}
\def\V{\bigvee} \def\U{\bigwedge}

\def\1{\oslash} \def\2{\oplus} \def\3{\otimes} \def\4{\ominus}
\def\5{\circ} \def\6{\odot} \def\7{\backslash} \def\8{\infty}
\def\9{\bigcap} \def\0{\bigcup} \def\+{\pm} \def\-{\mp}
\def\la{\langle} \def\ra{\rangle}

\def\Yl{Y_{h_{l}}^{\vv, i, M}} \def\Ylm{Y_{h_{l-1}}^{\vv, i, M}}
\def\tnp{t_{n+1}}
\def\L{\mathcal{L}}
\def\Plk{\mathcal{L}_{h_{l}}^{\vv, Y^k_n, M}}
\def\Plktwo{\mathcal{L}_{h_{l_2}}^{\vv, Y^k_n, M}}
\def\Plkm{\mathcal{L}_{h_{l-1}}^{\vv, Y^k_n, M}}
\def\Plj{\mathcal{L}_{h_{l}}^{\vv, Y^j_n, M}}
\def\Plr{\mathcal{L}_{h_{l}}^{\vv, Y^r_n, M}}
\def\Pljs{\mathcal{L}_{h_{l}}^{\vv, Y^{s,j}_n, M}}
\def\Pljm{\mathcal{L}_{h_{l-1}}^{\vv, Y^j_n, M}}
\def\Pl{\mathcal{L}_{h_{l}}^{\vv, Y_n, M}}
\def\Plm{\mathcal{L}_{h_{l-1}}^{\vv, Y_n, M}}
\def\Pls{\mathcal{L}_{h_{l}}^{\vv, Y^s_n, M}}
\def\Var{\rm Var} \def\Cov{\rm Cov}
\def\ff{\frac}


\def\tl{\tilde}
\def\trace{\hbox{\rm trace}}
\def\diag{\hbox{\rm diag}}
\def\for{\quad\hbox{for }}
\def\refer{\hangindent=0.3in\hangafter=1}

\newcommand\wD{\widehat{\D}}
\newcommand{\ka}{\kappa_{10}}

\title{
\bf Multilevel Monte Carlo  EM scheme for MV-SDEs with small noise\thanks{Supported }
 }

\author{
	{\bf  Ulises Botija-Munoz and Chenggui Yuan
	}\\
	\footnotesize{Department of Mathematics, Swansea University, Bay Campus, SA1 8EN, UK}\\
\footnotesize{Email: 942493@swansea.ac.uk, c.yuan@swansea.ac.uk}}

\date{}

\maketitle

\begin{abstract}
In this paper, we estimate the variance of two coupled paths  derived with the Multilevel Monte Carlo method combined with the Euler Maruyama discretization scheme for the simulation of McKean-Vlasov stochastic differential equations with small noise. The result often translates into a more efficient method than the standard Monte Carlo method combined with algorithms tailored to the small noise setting.   



\end{abstract}

\section{Introduction} An important problem in science is to approximate the value $\E[\Phi(X_T)]$ where $\{X_t\}_{0 \leq t \leq T}$ is the solution to an SDE and $\Psi:\RR^d \rightarrow \RR$. Among all the methods that allow us to compute the previous expectation, Monte Carlo simulation is arguably the more flexible. Its drawback is the high computational cost. Therefore a lot of effort has been placed to reduce this cost. In 2008, Giles, in a very relevant paper, \cite{g081}, proposed the multilevel Monte Carlo (MLMC) method which greatly reduces the computational cost of solving the problem $\E[\Phi(X_T)]$ with respect to the standard Monte Carlo (MC) method. If $\delta$ is the accuracy in terms of confidence intervals, the computation of $\E[\Phi(X_T)]$ where $X_T$ is simulated using the Euler-Maruyama (EM) method, has a computation cost (measured as the number of times that the random number generator is called) that scales as $\delta^{-3}$. In \cite{g081}, it is proved that the cost of the MLMC combined with the EM scheme scales like $\delta^{-2}(\log \delta)^2.$  Since then, numerous papers have appeared to customize, adapt and extend the principles of MLMC method to specific problems. One of these papers is \cite{ahs15}, where the authors applied the multilevel Monte Carlo framework to SDEs with small noise. They compare the computation cost derived from the standard MC method (combined with discretization algorithms tailored to the small noise setting) versus the multilevel Monte Carlo method combined with the Euler-Maruyama (EM) scheme.  They found that when $\delta \leq \e^2$, there is not benefit from using discretization methods customized for the small noise case. Moreover, if $\delta \geq e^{-\ff 1 \e}$, the EM scheme combined with the MLMC method leads to a cost $O(1).$ This is the same cost we would have with the standard MC method if we had $X_T$ as a formula of $W_T$, so no discretization method was required. In other words, the discretization method comes for free. Here, we extend the work from \cite{ahs15} to McKean-Vlasov SDEs (MV-SDEs) with small noise and we obtained the same estimate for the variance of two coupled paths. This means that the additional McKean-Vlasov component does not add computational complexity (per equation in the system of particles) and their conclusion about the computational cost of the method remains valid in our case. If we have a system of SDEs with $M$ particles, the total complexity is $M$ times the complexity of simulating one particle.
The MV-SDE with small noise that we will be working on in this paper, has the form
\begin{equation}\label{2.0}
	{\mbox d}X^{\vv}(t)=f(X^{\vv}(t),\mathcal{L}_{t}^{X}){\mbox d}t+ \vv g(X^{\vv}(t),\mathcal{L}_{t}^{X}){\mbox d}W(t),t \geq 0
\end{equation}
with initial data $X(0)=x_0$, where $\vv\in(0,1)$, $\mathcal{L}_{t}^{X}$ is the law (or distribution) of~$X(t)$, and
\begin{displaymath}
	f:\RR^d \times \mathcal P_2(\RR^{d})
	\rightarrow \RR^d  \mbox{ and }g: \RR^d \times \mathcal P_2(\RR^{d}) \rightarrow \RR^{d \times \bar d}.
\end{displaymath}

The pionering work on McKean-Vlasov SDEs is due to McKean on his work on the Boltzmann equation, \cite{mckeanII}. Since then, MV-SDEs have been used extensively in in biological systems, financial engineering and physics, \cite{Bala},\cite{Buck}, \cite{Gu}, \cite{Erban}. The existence and uniqueness theory for strong solutions of MV-SDEs with coefficients satisfying the Lipschitz condition is well-established, see, e.g., \cite{Szn}. Due to the propagation of chaos result \cite{mckeanI}, Equation \eqref{2.0} can be regarded as the limit of the following interacting particle system
\begin{align} \label{particleSystem}
	\mathrm{d}X^{\vv,i,M}(t)
	=f(X^{\vv, i,M}(t),\frac{1}{M}\sum\limits_{j=1}^{M}\delta_{X^{\vv, j,M}(t)})\mathrm{d}t
	+\vv g(X^{\vv,i,M}(t),\frac{1}{M}\sum\limits_{j=1}^{M}\delta_{X^{\vv, j,M}(t)})\mathrm{d}W^{i}(t),~~~~t\in[0,T].
\end{align}
Our main task in the rest of the paper is to discretize \eqref{particleSystem} using the EM scheme and estimate the variance of two coupled paths in the Multilevel Monte Carlo setting. This directly translate into the computational cost of solving $\E[\Phi((X^{\vv,i,M}(T)],$ see \cite{ahs15} for details.

\section{Preliminaries}

\subsection{Notation}

 Throughout this paper, unless otherwise specified, we let $(\Omega,
  {\cal{F}},\{{\cal{F}}_{t}\}_{t\ge 0}, \mathbb{P})$ be a complete probability space
  with a filtration
  $\{{\cal F}_t\}_{t\ge 0}$ satisfying
 the usual conditions.
 Let $W(t)=$ $(W_1(t),\ldots,W_{\bar d}(t))^T$ be an $\bar d$-dimensional Brownian motion defined on
 the probability space. For any $q>0$, let~$L^{q}=L^{q}(\Omega;\RR^{d})$ be the family of~$\RR^{d}$-valued random variables $Z$ with $\mathbb{E}[|Z|^q]<+\infty$. Let $\mathcal{L}^{Z}$ denote the probability law (or distribution) of a random variable~$Z$. $\delta_{x}(\cdot)$~denotes the Dirac delta measure concentrated at a point~$x\in\RR^{d}$.
 For $q\geq1$, we denote by $\mathcal{P}_q(\RR^{d})$ the set of probability measures on~$\RR^{d}$ with finite $q$th moments, and define 
 	\begin{equation}\label{qmoments}
 	W_q(\mu):=\left(\int_{\RR^{d}}|x|^q\mu(\mathrm{d}x)\right)^{\frac{1}q}, \quad \forall  \mu \in \mathcal{P}_q(\RR^{d}).
 	\end{equation}
 We assume that $(\Omega,{\cal{F}},\{{\cal{F}}_{t}\}_{t\ge 0}, \mathbb{P})$ is atomless so that, for any $\mu \in \mathcal P_2(\RR^d)$, there exists a random variable $X \in L^2(\Omega, \F, \PP; \RR^d)$ such that $\mu = \L^X$.
 Consider the MV-SDE with small noise \eqref{2.0}.
Let $f_i$ be the $i^{th}$ component of $f$. Then  for $x \in \RR^d$ and $\mu \in \mathcal P_2(\RR^d),$ we denote
\begin{align*}
\nabla f_i(x,\mu) &:=\left(\ff {\partial f_i(x,\mu)}{\partial x_1},...,\ff {\partial f_i(x,\mu)}{\partial x_d}\right),	\\
\nabla^2 f_i(x,\mu) &:=
\begin{bmatrix} 
	\ff {\partial^2 f_i(x,\mu)}{\partial x^2_1} & ... & \ff {\partial^2 f_i(x,\mu)}{\partial x_1 x_d} \\
	\vdots & \vdots & \vdots\\
	\ff {\partial^2 f_i(x,\mu)}{\partial x_d x_1} & ... & \ff {\partial^2 f_i(x,\mu)}{\partial x^2_d} \\
\end{bmatrix}.
\end{align*}

\begin{lemma}\label{defn2.1} \cite{carmona} (  Wasserstein~Distance
)
Let~$q\geq1$. Define
\begin{equation}\label{wasser}
\mathbb{W}_q(\mu,\nu):=\inf_{\pi \in \mathcal{D}(\mu,\nu)}\bigg\{\int_{\RR^{d}}|x-y|^q\pi(\mathrm{d}x,\mathrm{d}y)\bigg\}^{\frac{1}q},
~\mu,\nu\in\mathcal{P}_q(\RR^{d}),
\end{equation}
where~$\mathcal{D}(\mu,\nu)$~is the set of all couplings for~$\mu$~and~$\nu$. Then~$\mathbb{W}_q$~is a distance on~$\mathcal{P}_q(\RR^{d})$.
\end{lemma}
\begin{lemma}\label{l3.1}{\cite{carmona} }~~For any~$\mu\in\mathcal{P}_{2}(\RR^{d})$, $\mathbb{W}_{2}(\mu,\delta_{0})=W_{2}(\mu)$.
\end{lemma}

\subsection{Lions Derivatives}

In this subsection, we will give the definition of Lions derivative (or $L$-derivative) for a function $u:\mathcal{P}(\RR^{d}) \rightarrow \RR^{d}$  as introduced in \cite{car}. Given $(\Omega, \F, \PP),$ an atom is $A \in \F$ such that  $\PP(A)>0$ and for any $B \in \F, B \subset A, \PP(A)>\PP(B)$, we have that $\PP(B)=0.$	  
	
\begin{defn}
	We say that $u:\mathcal{P}(\RR^{d}) \rightarrow \RR^{d}$ is $L$-differentiable at $\mu \in \P(\RR^d)$ if there is an atomless, probability space $(\Omega, \F, \PP)$ and an $X \in L^2(\Omega, \F, \PP;\RR^d)$ such that $\mu = \mathcal L(X)$ and the function $U: L^2(\Omega, \F, \PP;\RR^d) \rightarrow \RR$ given by $U(X):=u(\L(X))$ is Frechet differentiable at $X$.
\end{defn}
We recall that $U$ is Frechet differentiable at $X$ means that there exists a continuous mapping $DU(X):L^2(\Omega, \F, \PP;\RR^d) \rightarrow \RR$ such that for any $Y \in L^2(\Omega, \F, \PP;\RR^d)$ 
$$U(X+Y)-U(X)=DU(X)(Y)+o(|Y|_{L^2}), \quad \text{as } |Y|_{L^2} \rightarrow 0.$$ 
Since $DU(X) \in  L^2(\Omega, \F, \PP;\RR^d)$, by Riesz representation theorem, there exists a $\PP$-a.s. unique variable $Z \in L^2(\Omega, \F, \PP;\RR^d)$ such that for any $Y \in L^2(\Omega, \F, \PP;\RR^d)$
$$DU(X)(Y) = \langle Y,Z \rangle_{L^2} = \E[YZ].$$ 
Cardaliaguet showed in \cite{car} that there exists a Borel measurable function $h:\RR^d \rightarrow \RR^d$ which only depends on the distribution $\L(X)$ rather than $X$ itself such that $Z=h(X)$. Thus, for $X \in L^2(\Omega, \F, \PP;\RR^d)$
$$u(\L(Y))-u(\L(X))=\E[h(X)(Y-X)]+o(|Y-X|_{L^2}).$$
We call $\partial_\mu u(\L(X))(y):=h(y), y \in \RR^d$ the $L$-derivative of $u$ at $\L(X),X \in L^2(\Omega, \F, \PP;\RR^d).$ 	 

Let $\bar u:\mathcal{P}(\RR^{d}) \rightarrow \RR$ be $L$-differentiable. Then by the mean value theorem (see chapter 5 in \cite{carmona}), for any two $d$-dimensional random variables $X$ and $X'$, there exists a $\theta \in [0,1]$ such that 
\begin{equation} \label{mvt}
	\bar u(\L(X))- \bar u(\L(X'))=\E[\langle \partial_\mu \bar u(\L(\theta X +(1-\theta)X'))(\theta X +(1-\theta)X'),(X-X')\rangle].
\end{equation}

\section{Multilevel Monte Carlo EM scheme for MV-SDEs with small noise}\label{sec2}
  We shall impose the following hypothesis on the functions $f$ and $g$:

\begin{assp}\label{a1} There exists a positive constant~ $K>0$ such that
\begin{align}\label{a1a}
|f(x,\mu)-f(y,\nu)|^2\vee |g(x,\mu)-g(y,\nu)|^2
&\leq K\big(|x-y|^{2}+\mathbb{W}^{2}_{2}(\mu,\nu)\big), 
\end{align}
hold for any $x,~y\in\RR^{d}$,~$\mu,~\nu\in\mathcal{P}_{2}(\RR^{d})$.
Furthermore there exists a positive constant $\bar K$ such that
$$|\nabla f(x,\mu)|^2 \vee |\nabla^2 f(x,\mu)|^2 \vee |\partial_{\mu} f(x,\mu)(y)|^2 \vee |\partial^2_{\mu} f(x,\mu)(y)|^2 \leq \bar K$$
for all $x,~y\in\RR^{d}$,~$\mu \in\mathcal{P}_{2}(\RR^{d})$.
In addition, there exists a positive constant $K$ such that
\begin{equation}\label{a2a}
|\partial_{\mu} f(x,\mu)(y)-\partial_{\mu} f(\bar x,\nu)(\bar y)|^2 \leq K\big(|x-\bar x|^2 +|y-\bar y|^2 +\mathbb W^2_2(\mu,\nu)\big).
\end{equation}
for all $x,y,\bar x, \bar y \in\RR^{d}$,~$\mu,\nu \in\mathcal{P}_{2}(\RR^{d}).$
\end{assp}

\begin{lemma}
Let Assumption \ref{a1} hold. Then, for any $T>0$ and $p\ge 2,$ we have
$$
\E \left[\sup_{0\le t\le T}|X^\vv(t)|^p\right] \le C.
$$
\end{lemma}
The proof of this lemma is standard, we omit it here.

\begin{remark}\label{onem}
{\rm Assumption \ref{a1}  implies the existence and uniqueness of equation \eqref{2.0}. Moreover, if Assumption \ref{a1}, then 
\begin{displaymath}
|f(x,\mu)|^2 \vee |g(x, \mu)|^2\le \beta(1+|x|^2 + W^{2}_{2}(\mu)),
\end{displaymath}
where $\beta=2\max\{1, |f(0,\delta_{0})|,|g(0,\delta_{0})|\}$, and for any $x \in \RR^d$ and $\mu\in\mathcal{P}_{2}(\RR^{d}).$
\begin{displaymath}
\langle x-y,f(x,\mu)-f(y,\nu)\rangle  \le \alpha(|x-\bar{x}|^2+\mathbb{W}^{2}_{2}(\mu,\nu)),
\end{displaymath}
where $\bar{\alpha}=\frac{1}{2}(1+K)$.
}
\end{remark}

\subsection{Stochastic Particle Method}\label{s3.2}

In this subsection, we make use of the stochastic particle method \cite{bossy} to approximate the MV-SDDE \eqref{2.0}. For any~$ i\in\mathbb{N}$,
$\{W^{i}(t)\}_{t\in[0,T]}$~is~an $\bar d$-dimension Brownian motion.  Assume $\{W^{1}(t)\}, \{W^{2}(t)\}, \cdots$~are  independent and $x^{1}, x^{2}, \cdots$ are independent and identically distributed ($i.i.d.$).  Let~$\{X^{\e,i}(t)\}_{t\in[0,T]}$~be the unique solution to MV-SDE
\begin{equation} \label{eq3.32}
\mathrm{d}X^{\vv, i}(t)=f(X^{\vv, i}(t),\mathcal{L}_{t}^{X^{\vv, i}})
\mathrm{d}t+\vv g(X^{\vv, i}(t),\mathcal{L}_{t}^{X^{\vv, i}})W^{i}(t),
\end{equation}
with the initial condition~$X^{i}_0 =x^{i}$
 and~$\mathcal{L}_{t}^{X^{\vv, i}}$~is the law of~$ X^{\vv, i}(t)$.  One can see  that~$ X^{\vv, 1}(t), X^{\vv, 2}(t), \cdots$ are i.i.d.for $t\geq 0$.

For any  $M\in\mathbb{N},~1\leq i\leq M$, let~$ X^{\e,i,M}(t) $~be the solution of SDE
\begin{align} \label{eq3.33}
\mathrm{d}X^{\vv,i,M}(t)
=f(X^{\vv, i,M}(t),\mathcal{L}_{t}^{\vv,X,M})\mathrm{d}t
+\vv g(X^{\vv,i,M}(t),\mathcal{L}_{t}^{\vv,X,M})\mathrm{d}W^{i}(t),~~~~t\in[0,T],
\end{align}
with the initial condition~$ X^{\e,i,M}_0 =x^{i} $, where 
$\mathcal{L}_{t}^{\vv,X,M}:=\frac{1}{M}\sum\limits_{j=1}^{M}\delta_{X^{\vv, j,M}(t)}$.
We prepare a path-wise propagation of chaos result on SDEs~$(\ref{eq3.33})$.

\begin{lemma} \cite{li} \label{le3.6}
If Assumption~$\ref{a1}$~holds, then
\begin{align*}
\displaystyle\sup_{1\leq i\leq M}\mathbb{E}\big[\sup_{0\leq t\leq T}|X^{\e,i}(t)-X^{\e,i,M}(t)|^{2}\big]\leq C\left\{
\begin{array}{lll}
M^{-\frac{1}{2}},~~~&1\leq d<4,\\
M^{-\frac{1}{2}}\log(M),~~~&d=4,\\
M^{-\frac{d}{2}},~~~&4<d,
\end{array}
\right.
\end{align*}
where $C$ is independent of $M$.
\end{lemma}

\subsection{The  EM Scheme for MV-SDEs with small noise} \label{sec2.1}
We now introduce the EM scheme for \eqref{2.0}. Given any time $T>0$, assume that  there exists a positive integer such that $h=\frac{T}{m}$, where $h\in (0,1)$ is the step size.  Let  $t_n=nh$ for $n \ge 0$. Compute the discrete approximations $Y_{h, n}^{\vv, i, M}=Y_{h}^{\vv, i, M}(t_n)$ by setting $Y_h^{\vv, i, M}(0)=x_0$ and  forming
\begin{equation}\label{discrete}
\begin{split}
Y_{h, n+1}^{\vv, i, M}=Y_{h, n}^{\vv, i, M}
+ f(Y_{h, n}^{\vv, i, M}, \mathcal{L}_h^{\vv, Y_n,M} )h+\vv g(Y_{h, n}^{\vv, i, M}, \mathcal{L}_h^{\vv, Y_n,M} )\Delta W^i(t_n),
\end{split}
\end{equation}
where    $\mathcal{L}_h^{\vv, Y_n,M} =\frac{1}{M}\sum\limits_{j=1}^{M}\delta_{Y^{\vv, j,M}_{h, n}}$ and $\Delta W(t_n)=W(t_{n+1})-W(t_n)$.

Let
\begin{align}\label{eq4.3}
Y_h^{\vv, i,M}(t)=Y^{\vv, i,M}_{h,k}, ~~~~ t\in[t_k, t_{k+1}).
\end{align}
For convenience, we define $\mathcal{L}^{\vv, Y,M}_{h, t}=\frac{1}{M}\sum\limits_{j=1}^{M}\delta_{Y_h^{\vv, i,M}(t)}$ and $\eta_h(t):=\lfloor t/h\rfloor h$~for~$t\geq 0$. Then one observes~$\mathcal{L}^{\vv, Y,M}_{h, t}=\mathcal{L}^{\vv, Y,M}_{h, \eta_h(t)}=\mathcal{L}^{\vv, Y_k,M}_{h}$, for $t\in [t_k, t_{k+1})$.
We now define the EM continuous approximate solution as follows:
\begin{equation}\label{EMsol}
\bar{Y}_h^{\vv, i,M}(t)=x^{i}+\int^{t}_{0}f(Y_h^{\vv, i,M}(s),\mathcal{L}_{h, s}^{\vv, Y,M})\mathrm{d}s
+\e\int^{t}_{0}g (Y_h^{\vv, i,M}(s),\mathcal{L}_{h,s}^{\vv, Y,M})\mathrm{d}W^{i}(s), ~ t\ge 0.
\end{equation}

\begin{lemma}\label{0pmoment}
Let Assumption \ref{a1} hold. Then, for any $T>0$ and $p\ge 2,$ we have
$$\E \left[\sup_{0\le t\le T}|\bar Y_{h}^{\vv, i,M}(t)|^p\right] \le C.$$
\end{lemma}
The proof of this lemma is standard, we omit it here.

\begin{lemma}\label{womiss}
Let Assumption \ref{a1} hold. Then, for any $p\ge2,$ we have
\begin{equation*}
\sup_{0\le t \le T}\E [|\bar Y^{\vv, i,M}_{h}(t)- Y_{h}^{\vv, i,M}(t)|^p]\le Ch^p+  C \vv^ph^{p/2}.
\end{equation*}
\end{lemma}
{\bf Proof.} Let $n$ be such that $t_n \leq t \leq t_{n +1}.$ From \eqref{EMsol} we have
$$\bar{Y}_h^{\vv, i,M}(t)- Y_h^{\vv, i,M}(t)=\int^{t}_{t_n}f(Y_h^{\vv, i,M}(s),\mathcal{L}_{h, s}^{\vv, Y,M})ds
+\e\int^{t}_{t_n}g (Y_h^{\vv, i,M}(s),\mathcal{L}_{h,s}^{\vv, Y,M})dW^{i}(s)$$
By Remark \ref{onem} and the BDG inequality, one has
\begin{equation*}
\begin{split}
\E|\bar Y^{\vv, i,M}_{h}(t)- Y_{h}^{\vv, i,M}(t)|^p\le&2^{p-1}h^{p-1}\E\int_{t_n}^t|f(Y_h^{\vv, i,M}(s),\mathcal{L}_{h, s}^{\vv, Y,M})|^p ds\\
&+\vv^ph^{\frac{p}{2}-1}\E\int_{t_n}^t|g(Y_h^{\vv, i,M}(s),\mathcal{L}_{h, s}^{\vv, Y,M})|^p ds\\
\le& Ch^p+ C\vv^ph^{\frac{p}{2}}.
\end{split}
\end{equation*}
The proof is therefore complete. \hfill $\Box$

We now reveal the error between the numerical solution \eqref{EMsol} and the exact solution \eqref{2.0}.

\begin{theorem}\label{th0}
Let Assumption \ref{a1} hold, assume that $\Psi:\mathbb{R}^d\rightarrow\mathbb{R}$ has continuous second order derivative and there exists a constant $C$ such that
\begin{equation*}
\begin{split}
\left|\frac{\partial \Psi}{\partial x_i}\right|\le C
\end{split}
\end{equation*}
for any $i=1,2,\cdots,d$. Then we have
\begin{equation*}
\begin{split}
\sup\limits_{0\le t\le T}\mathbb{E}|\Psi(X^{\vv, i, M}(t))-\Psi(\bar Y_h^{\vv, i, M}(t))|^2
=Ch^2+  C h\e^2.
\end{split}
\end{equation*}
\end{theorem}
{\bf Proof. }  By Assumption \ref{a1} and Lemma \ref{womiss}, one can see that
\begin{equation}\label{new}
\begin{split}
&\sup\limits_{0\le t\le T}\mathbb{E}|X^{\vv, i, M}(t)-\bar Y_h^{\vv, i, M}(t)|^2\\
&\le 2T\mathbb{E}\int_0^T|f(X^{\vv, i,M}(s),\mathcal{L}_{s}^{\vv,X,M})-f(Y_h^{\vv, i, M}(s), \mathcal{L}_{h, s}^{\vv, Y,M})|^2 ds\\
&+8\sqrt{T}\vv^2\mathbb{E}\int_0^T|g(X^{\vv, i,M}(s),\mathcal{L}_{s}^{\vv,X,M})-g(Y_h^{\vv, i, M}(s), \mathcal{L}_{h, s}^{\vv, Y,M})|^2 ds\\
&\le  2TK\mathbb{E}\int_0^T\big(|X^{\vv, i,M}(s)-Y_h^{\vv, i, M}(s)|^2+W_2^2(\mathcal{L}_{s}^{\vv,X,M}, \mathcal{L}_{h, s}^{\vv, Y,M}\big)  \\
&+ 8K\sqrt{T}\vv^2\mathbb{E}\int_0^T\big(|X^{\vv, i,M}(s)-Y_h^{\vv, i, M}(s)|^2+W_2^2(\mathcal{L}_{s}^{\vv,X,M}, \mathcal{L}_{h, s}^{\vv, Y,M}\big)\\
&\le 4TK\mathbb{E}\int_0^T|X^{\vv, i,M}(s)- \bar Y_h^{\vv, i, M}(s)|^2 ds + 4TK\mathbb{E}\int_0^T|\bar Y^{\vv, i,M}(s)- Y_h^{\vv, i, M}(s)|^2 ds \\ 
&+16K\sqrt{T}\vv^2\mathbb{E}\int_0^T |X^{\vv, i,M}(s)- \bar Y_h^{\vv, i, M}(s)|^2 ds +16K\sqrt{T}\vv^2\mathbb{E}\int_0^T |\bar Y^{\vv, i,M}(s)-Y_h^{\vv, i, M}(s)|^2 ds\\
&\le Ch^2+ C\vv^2h+C\vv^2\int_0^T\sup\limits_{0\le t\le s}\mathbb{E}|X^{\vv, i,M}(s)-\bar Y_h^{\vv, i, M}(s)|^2 ds+C\vv^2 h^2+C\vv^4h.
\end{split}
\end{equation}
The Gronwall inequality implies that
\begin{equation*}
\begin{split}
\sup\limits_{0\le t\le T}\mathbb{E}|X^{\vv, i,M}(t)-\bar Y_h^{\vv, i, M}(t)|^2\le Ch^2+\vv^2h.
\end{split}
\end{equation*}
Since $\Psi$ has continuous bounded first order derivative, we immediately get
\begin{equation*}
\begin{split}
\sup\limits_{0\le t\le T}\mathbb{E}|\Psi(X^{\vv, i,M}(t))-\Psi(\bar Y_h^{\vv, i, M}(t))|^2\le C\sup\limits_{0\le t\le T}\mathbb{E}|X^{\vv, i,M}(t)-\bar Y_h^{\vv, i, M}(t)|^2.
\end{split}
\end{equation*}
The desired result then follows. \hfill $\Box$

In the next corollary, we are going to use different stepsize to define the numerical solutions.
\begin{coro}
Assume that the conditions of Theorem \ref{th0} hold. Let $M\ge 2, l\ge 1$, $h_l=T\cdot M^{-l}, h_{l-1}=T\cdot M^{-(l-1)}$.  Then
\begin{equation*}
\begin{split}
\max_{0\le n<M^{l-1}}{\rm Var}(\Psi(\bar Y_{h_l}^{\vv, i, M}(t_n))-\Psi(\bar Y_{h_{l-1}}^{\vv, i, M}(t_n))\le Ch_{l-1}^2+C\vv^2h_{l-1}.
\end{split}
\end{equation*}
\end{coro}
{\bf Proof.} For $0\le n\le M^{l-1}-1$, by Theorem \ref{th0},
\begin{equation*}
\begin{split}
&{\rm Var}(\Psi(\bar Y_{h_l}^{\vv, i, M}(t_n))-\Psi(\bar Y_{h_{l-1}}^{\vv, i, M}(t_n))\le2 \mathbb{E}|\Psi(\bar Y_{h_l}^{\vv, i, M}(t_n))-\Psi(\bar Y_{h_{l-1}}^{\vv, i, M}(t_n))|^2\\
\le&4\mathbb{E}|\Psi(\bar Y_{h_l}^{\vv, i, M}(t_n))-\Psi(X^{\vv, i,M}(t))|^2+2\mathbb{E}|\Psi(X^{\vv, i,M}(t))-\Psi(\bar Y_{h_{l-1}}^{\vv, i, M}(t_n))|^2\\
\le&Ch_{l-1}^2+C\vv^2h_{l-1}.
\end{split}
\end{equation*}
\hfill $\Box$

The following lemma is presented here because it applies to any EM scheme but it will only be use later when estimating the variance of coupled processes in the Multilevel Monte Carlo setting.

Define $\eta_h(s):=\lfloor s/h \rfloor$ where $\lfloor \cdot \rfloor$ is the integer-part function. Let $z_h$ be the deterministic solution to
\begin{equation}\label{z_h}
z_h(t) = X(0) + \int_0^t f(z_h(\eta_h(s)),\delta_{z_h(s)})ds,
\end{equation}
which is the Euler approximation to the ODE obtained from \eqref{2.0} when $\e$ is set to zero.
\begin{lemma}\label{z}
	For any $T>0$ we have
	\begin{equation}\label{ODE}
	\E[\sup_{0 \leq s \leq T} |\bar Y_h^{\vv, i, M}(s) - z_h(s)|^2] \leq C\e^2.
	\end{equation}
\end{lemma}	
{\bf Proof.} Using \eqref{EMsol} and \eqref{ODE}, using the fact that $|a+b|^2 \leq a^2 + b^2$ and the Cauchy-Schwarz inequality we have that for every $t \leq T$
\begin{align*}
	|&\bar{Y}_h^{\vv, i,M}(t)- z_h(t)|^2\\
	&=\left|\int^t_0(f(Y_h^{\vv, i,M}(s),\mathcal{L}_{h, s}^{\vv, Y,M}) -f(z_h(\eta_h(s)),\delta_{z_h(s)})))ds
	+\e\int^{t}_{0}g (Y_h^{\vv, i,M}(s),\mathcal{L}_{h,s}^{\vv, Y,M})\mathrm{d}W^{i}(s)\right|^2 \\
	&\leq 2T \int^t_0|f(Y_h^{\vv, i,M}(s),\mathcal{L}_{h, s}^{\vv, Y,M}) -f(z_h(\eta_h(s)),\delta_{z_h(s)}))|^2 ds + 2\e^2\left|\int^{t}_{0}g (Y_h^{\vv, i,M}(s),\mathcal{L}_{h,s}^{\vv, Y,M})\mathrm{d}W^{i}(s)\right|^2.
\end{align*}
By the BDG inequality we have that
$$\E \left[\sup_{0 \leq s \leq t}\left|\int^{t}_{0}g (Y_h^{\vv, i,M}(s),\mathcal{L}_{h,s}^{\vv, Y,M})\mathrm{d}W^{i}(s)\right|^2 \right] \leq 4 \int_0^t \E[|g (Y_h^{\vv, i,M}(s),\mathcal{L}_{h,s}^{\vv, Y,M})|^2]ds.$$
Thus by Assumption \ref{a1} one can see that
\begin{align*}
	\E&[\sup_{0 \leq s \leq t}|\bar{Y}_h^{\vv, i,M}(t)- z_h(t)|^2 ] \leq 2TK \int^t_0(\E[\sup_{0 \leq s \leq r}|\bar{Y}_h^{\vv, i,M}(s)- z_h(s)|^2]+\sup_{0 \leq s \leq r}\mathbb W^{2}_{2}(\mathcal{L}_{h, s}^{\vv, Y,M},\delta_{z_h(s)})) dr\\
	  &+ 8T\e^2 \beta \int^{t}_{0}\E[(1+|\bar{Y}_h^{\vv, i,M}(s)|^2 + W^{2}_{2}(\mathcal{L}_{h, s}^{\vv, Y,M})] ds.
\end{align*}
Using \eqref{qmoments}, \eqref{wasser} and Lemma \ref{0pmoment} we have that for all $0 \leq t \leq T$
\begin{align*}
	\E[\sup_{0 \leq s \leq t}|\bar{Y}_h^{\vv, i,M}(t)- z_h(t)|^2 ] \leq C \e^2 +  C \int^t_0\E[\sup_{0 \leq s \leq r}|\bar{Y}_h^{\vv, i,M}(s)- z_h(s)|^{2}] dr. 
\end{align*}
The final result is obtained by applying the Gronwall inequality.  
\hfill $\Box$

\subsection{The Multilevel Monte Carlo  EM Scheme}\label{sec2.2}
We now define the multilevel Monte Carlo EM scheme. Given any $T>0$, let $N\ge 2, l \in \{0,...,L\}$, where $L$ is a positive integer that will be determined later. Let $h_l=T\cdot N^{-l}, h_{l-1}=T\cdot N^{-(l-1)}$. 

For step sizes $h_{l}$ and $h_{l-1}$ the EM continuous approximate solutions are respectively
\begin{equation}\label{c1}
	\begin{split}
		\bar{Y}_{h_l}^{\vv, i,M}(t)=x^{i}+\int^{t}_{0}f(Y_{h_l}^{\vv, i,M}(s),\mathcal{L}_{{h_l}, s}^{\vv, Y,M})\mathrm{d}s
		+\int^{t}_{0}g (Y_{h_l}^{\vv, i,M}(s),\mathcal{L}_{{h_l},s}^{\vv, Y,M})\mathrm{d}W^{i}(s), 
	\end{split}
\end{equation}
and 
\begin{equation}\label{c2}
	\begin{split}
		\bar{Y}_{h_{l-1}}^{\vv, i,M}(t)=x^{i}+\int^{t}_{0}f(Y_{h_{l-1}}^{\vv, i,M}(s),\mathcal{L}_{{h_{l-1}}, s}^{\vv, Y,M})\mathrm{d}s
		+\int^{t}_{0}g (Y_{h_{l-1}}^{\vv, i,M}(s),\mathcal{L}_{{h_{l-1}},s}^{\vv, Y,M})\mathrm{d}W^{i}(s).
	\end{split}
\end{equation}
We now construct the discrete version of the previous approximate solutions using the same Brownian motion for both processes. We say that the two processes are coupled.   
For $n\in \{0, 1, \ldots, N^{l-1}-1\}$ and $k\in \{0, \ldots, N\}$, let
$$
t_n=nh_{l-1} \mbox{ and } t_n^k=nh_{l-1}+kh_l.
$$
This means we divide the interval $[t_n, t_{n+1}]$ into $N$ equal parts by ${h_l}$ with $t_n^0=t_n, t_n^{N}=t_{n+1}.$ 
For  $n\in \{0, 1, \ldots, N^{l-1}-1\}$ and $k\in \{0, \ldots, N-1\}$, let
\begin{equation}\label{MMC EM scheme_1}
\begin{split}
Y_{h_{l}}^{\vv, i, M}(t_n^{k+1})=Y_{h_l}^{\vv, i, M}(t_n^{k})
+ f(Y_{h_l}^{\vv, i, M}(t_n^{k}), \mathcal{L}_{h_l}^{\vv, Y_n^{k},M} )h_l+\vv \sqrt{h_l}g(Y_{h_l}^{\vv, i, M}(t_n^{k}), \mathcal{L}_{h_l}^{\vv, Y_n^{k},M})\Delta \xi_n^k,
\end{split}
\end{equation}
where $\mathcal{L}_{h_l}^{\vv, Y_n^{k},M}=\frac{1}{M}\sum_{j=1}^M\delta_{Y_{h_l}^{\vv, j, M}(t_n^{k})}$, the random vector $\Delta \xi_n^k\in\mathbb{R}^{\bar d}$ has independent components, and each component is distributed as $\mathcal N(0, 1).$ Therefore, to simulate  $Y_{h_{l}}^{\vv, i, M},$ we use
\begin{equation}\label{MMC EM scheme_2}
\begin{split}
Y_{h_{l}}^{\vv, i, M}(t_{n+1})=Y_{h_l}^{\vv, i, M}(t_n)+ \sum_{k=0}^{N-1}f(Y_{h_l}^{\vv, i, M}(t_n^{k}), \mathcal{L}_{h_l}^{\vv, Y_n^{k},M} )h_l+\vv \sqrt{h_l}\sum_{k=0}^{N-1}g(Y_{h_l}^{\vv, i, M}(t_n^{k}), \mathcal{L}_{h_l}^{\vv, Y_n^{k},M})\Delta \xi_n^k.
\end{split}
\end{equation}
To simulate $Y_{h_{l-1}}^{\vv, i, M},$ we use
\begin{equation}\label{MMC EM scheme_3}
\begin{split}
Y_{h_{l-1}}^{\vv, i, M}(t_{n+1})
&=Y_{h_{l-1}}^{\vv, i, M}(t_{n})+f(Y_{h_{l-1}}^{\vv, i, M}(t_{n}), \mathcal{L}_{h_{l-1}}^{\vv, Y_n,M} )h_{l-1}\\
&+\vv \sqrt{h_l}g(Y_{h_{l-1}}^{\vv, i, M}(t_{n}), \mathcal{L}_{h_{l-1}}^{\vv, Y_n, M})\sum_{k=0}^{N-1}\Delta \xi_n^k,
\end{split}
\end{equation}
where $\mathcal{L}_{h_{l-1}}^{\vv, Y_n,M}=\frac{1}{M}\sum_{j=1}^M\delta_{Y_{h_{l-1}}^{\vv, j, M}(t_n)}$.

The following theorem is the main result of this section.

\begin{theorem}\label{2ndMom2Paths} 
Let Assumption \ref{a1} hold. Then it holds that
\begin{equation*}
\max_{0\le n<M^{l-1}}\E [|Y_{h_{l}}^{\vv, i, M}(t_{n})-Y_{h_{l-1}}^{\vv, i, M}(t_{n})|^2]\le C N^2 h_l^2+ \bar C\vv^4 N h_l.
\end{equation*}
\end{theorem}

In order to prove Theorem \ref{2ndMom2Paths}, we need a few lemmas.

\begin{lemma} \label{lemma kn}
Let $0<p\leq 4$. Then 
$$\max_{\substack{0 \leq n \leq N^{l-1} \\ 1 \leq k \leq N}} \E[|\Yl(t^k_n) - \Yl(t_n)|^p] \leq C_1N^p h_l^p +C_2 N^{p/2}   h_l^{p/2} \e^p,   $$
where $C$ and $C$ are positive constants that only depend on $\beta, T,m$ and $X^{\e}(0)$ ($\beta$ from Remark \ref{onem}). 
\end{lemma}

{\bf Proof.} Let $p = 4$. From \eqref{MMC EM scheme_1} we  have that
\begin{equation}\label{Y^k-Y}	
\Yl(t^k_n) - \Yl(t_n) = \sum_{j=0}^{k-1}f(\Yl(t_n^j),\Plj)h_l + \e \sqrt{h_l} \sum_{j=0}^{k-1}g(\Yl(t_n^j),\Plj)\D \xi_n^j. 
\end{equation}
Hence, we obtain
\begin{align}\label{eq1-lemma1}
	\E[|\Yl(t^k_n) &- \Yl(t_n)|^4] \\
&\leq 8 \E \left|\sum_{j=0}^{k-1}f(\Yl(t_n^j),\Plj)h_l\right|^4 + 8\E\left|\e \sqrt{h_l} \sum_{j=0}^{k-1}g(\Yl(t_n^j),\Plj)\D \xi_n^j \right|^4. \nonumber 
\end{align}
By Remark \ref{onem} and Lemma \ref{0pmoment} one can see that
\begin{align}\label{eq2-lemma1}
	\E &\left| \sum_{j=0}^{k-1} f(\Yl(t_n^j),\Plj)h_l\right|^4 \leq N^3 \sum_{j=0}^{k-1} \E \left| f(\Yl(t_n^j),\Plj)h_l \right|^4  \nonumber \\
	&\leq N^3 \sum_{j=0}^{k-1} \E\left[\left(\beta \left(1+|\Yl(t^j_n)|^2 + W_2^2(\Plj) \right)\right)^2\right] \nonumber \\
	&\leq 3 N^3 h_l^4 \beta^2 \sum_{j=0}^{k-1} \left(1+2\E|\Yl(t^j_n)|^4 \right) \leq C N^4 h_l^4 .  
\end{align}
Using the BDG inequality, Remark \ref{onem} and Lemma \ref{0pmoment}, we obtain
\begin{align}\label{eq3-lemma1}
	\E &\left|\e \sqrt{h_l} \sum_{j=0}^{k-1} g(\Yl(t_n^j),\Plj)\D \xi_n^j \right|^4 \leq C \e^4 \E \left[ \left| \sum_{j=0}^{k-1} | g(\Yl(t_n^j),\Plj)|^2 h_l \right|^2 \right] \nonumber  \\
	&\leq C \e^4 N h_l^2 \E \left[ \sum_{j=0}^{k-1} (| g(\Yl(t_n^j),\Plj)|^2)^2   \right] \nonumber \\
	&\leq C \e^4 N h_l^2  \sum_{j=0}^{k-1} \E\left[\left(\beta \left(1+|\Yl(t^j_n)|^2 + W_2^2(\Plj) \right)\right)^2\right] \leq  C   N^2 h_l^2 \e^4 \nonumber\\
\end{align}
The result for $p=4$ follows from substituting \eqref{eq2-lemma1} and \eqref{eq3-lemma1} into \eqref{eq1-lemma1}. For $0 < p < 4$, the result follows from Jensen's inequality.
\hfill $\Box$

\begin{lemma}\label{ABE}
	Let $f_{m}$ be the $m^{th}$ component of $f$. Then there exist $s,r \in [0,1]$ such that 
	$$f(\Yl(t_n^k),\Plk)-f(\Yl(t_n),\Plk)=A_k +  B_k +  E_k,$$
	where
\begin{align*}
	A_k &= (A_k^1,...,A_k^d)', B_k = (B_k^1,...,B_k^d)', E_k = (E_k^1,...,E_k^d)' \\
	A^m_k &:= \langle \nabla f_m(s \Yl(t_n^k) + (1-s)\Yl(t_n)),\Plk), h_l \sum_{j=0}^{k-1} f(\Yl(t^j_n),\Plj) \rangle, \\
	B^m_k &:= \langle \nabla f_m(\Yl(t_n),\Plk), \e \sqrt{h_l} \sum_{j=0}^{k-1} g(\Yl(t^j_n),\Plj) \D \xi_n^j  \rangle, \\
	E^m_k &:=  \langle \nabla^2 f_m(rs(\Yl(t^k_n)- \Yl(t_n))+ \Yl(t_n),\Plk)(\Yl(t^k_n)- \Yl(t_n))s, \\
	 &\e \sqrt{h_l} \sum_{j=0}^{k-1} g(\Yl(t^j_n),\Plj) \D \xi_n^j \rangle , m \in \{1,...,d\}.
\end{align*}
\end{lemma}
{\bf Proof.} By the mean value theorem there exists a $s \in [0,1]$ such that 
\begin{align*} 
f_m(\Yl &(t_n^k),\Plk)-f_m(\Yl(t_n),\Plk)  \\
&= \langle \nabla f_m(s \Yl(t_n^k) + (1-s)\Yl(t_n)),\Plk), (\Yl(t_n^k) - \Yl(t_n)) \rangle.  \nonumber
\end{align*}
Substituting \eqref{Y^k-Y} in the equation above yields
\begin{align}\label{eq1-ABE}
	f_m(\Yl &(t_n^k),\Plk)-f_m(\Yl(t_n),\Plk) \\
	&=\langle  \nabla f_m(s \Yl(t_n^k) + (1-s)\Yl(t_n)),\Plk),  \sum_{j=0}^{k-1}f(\Yl(t_n^j),\Plj)h_l \rangle \nonumber \\ 
	+& \langle \nabla f_m(s \Yl(t_n^k) + (1-s)\Yl(t_n)),\Plk), \e \sqrt{h_l} \sum_{j=0}^{k-1}g(\Yl(t_n^j),\Plj)\D \xi_n^j \nonumber \rangle.
\end{align}
Let $\nabla_q f_m$ denote the $q^{th}$ component of the vector function $\nabla f_m$. Applying the mean value theorem again with $y=s \Yl(t_n^k) +(1-s)\Yl(t_n),x=\Yl(t_n)$ and $g(z)=\nabla_q f_m(z, \Plk)$  ensures that there exists a $r \in [0,1]$ such that
\begin{align*}
\nabla_q f_m &(s \Yl(t_n^k) + (1-s)\Yl(t_n)),\Plk) = \nabla_q f_m(\Yl(t_n),\Plk) \\
&+  \langle \nabla(\nabla_q f_m)(rs (\Yl(t_n^k)-\Yl(t_n)) + \Yl(t_n),\Plk),(\Yl(t_n^k)-\Yl(t_n))s \rangle.
\end{align*}
Thus
\begin{align*}
	\nabla f_m &(s \Yl(t_n^k) + (1-s)\Yl(t_n)),\Plk) = \nabla f_m(\Yl(t_n),\Plk) \\
	&+   \nabla^2 f_m(rs (\Yl(t_n^k)-\Yl(t_n)) + \Yl(t_n),\Plk)(\Yl(t_n^k)-\Yl(t_n))s .
\end{align*}
Substituting the last equation into the second summand of the RHS of  \eqref{eq1-ABE} completes the proof.
\hfill $\Box$	

\begin{lemma}\label{bar ABE}
 There exist random variables $s,r:\Omega \rightarrow [0,1]$ such that
	$$f(\Yl(t_n),\Plk)-f(\Yl(t_n),\Pl) = \bar A_k + \bar E_k,$$
	where 
\begin{align*}
	\bar A_k &= (\bar A_k^1,..., \bar A_k^d)', \bar E_k = (\bar E_k^1,...,\bar E_k^d)' \\
	\bar A^m_k &:=\E[\langle \partial_{\mu} f_m(Z,\Pls)(Y_n^s), h_l \sum_{j=0}^{k-1} f(\Yl(t^j_n),\Plj) \rangle]_{Z=\Yl(t_n)} \\
	\bar E^m_k &:=\E[\langle \partial^2_{\mu} f_m(Z,\mathcal L_{h_l}^{\e,Y^{s,r}_n,M})(Y^{s,r}_n)(((\Yl(t^k_n)-\Yl(t_n))s, \\
	& \quad \quad \quad \quad \quad \quad \quad \quad \quad \quad \quad \quad \quad \e \sqrt{h_l} \sum_{j=0}^{k-1} g(\Yl(t^j_n),\Plj) \D \xi_n^j \rangle]_{Z=\Yl(t_n)}, \\
	Y_n^s &:= s \Yl(t_n^k) + (1-s)\Yl(t_n), \\
	Y^{s,r}_n &:= sr(\Yl(t^k_n)- \Yl(t_n))+ \Yl(t_n).
\end{align*} 
\end{lemma}
{\bf Proof.} Let $f_{m}$ be the $m^{th}$ component of $f$. A direct application of Equation \eqref{mvt} with $X = \Yl(t_n^k),X'=\Yl(t_n)$ and $\bar u(\mathcal L(\xi))=f_m(\Yl, \mathcal L^{\e,\xi,M}_{h_l})$  implies that there exists a random variable $s:\Omega \rightarrow [0,1]$ such that 
\begin{align}\label{eq1-bar ABE}
f_m(\Yl(t_n),&\Plk)-f_m(\Yl(t_n),\Pl) \\
&= \E[\langle \partial_{\mu} f_m(Z,\Pls)(Y_n^s),(\Yl(t^k_n)-\Yl(t_n)) \rangle]_{Z=\Yl(t_n)} \nonumber \\
&=\E[\langle \partial_{\mu} f_m(Z,\Pls)(Y_n^s),\sum_{j=0}^{k-1}f(\Yl(t_n^j),\Plj)h_l \rangle]_{Z=\Yl(t_n)} \nonumber \\
&+ \E[\langle \partial_{\mu} f_m(Z,\Pls)(Y_n^s), \e \sqrt{h_l} \sum_{j=0}^{k-1}g(\Yl(t_n^j),\Plj)\D \xi_n^j \rangle]_{Z=\Yl(t_n)}. \nonumber
\end{align} 
Let $\partial_{\mu,q} f_m$ be the $q^{th}$ component of the vector function $\partial_{\mu} f_m.$ Applying Equation \eqref{mvt} again with $X= s \Yl(t_n^k) + (1-s)\Yl(t_n)=:Y_n^s,X'=\Yl(t_n)$ and $\bar u(\mathcal L(\xi))= \partial_{\mu,q} f_m(\Yl(t_n),\mathcal L^{\e,\xi,M}_{h_l})(\xi)$, we find that there exists a random variable $r:\Omega \rightarrow [0,1]$ such that
\begin{align*}
	\partial_{\mu,q} &f_m(Z,\Pls)(Y_n^s)=\partial_{\mu,q} f_m(\Yl(t_n),\Pl)(\Yl(t_n)) \\
	&+\E[\langle\partial_{\mu}(\partial_{\mu,q}f_m)(Z,\mathcal L_{h_l}^{\e,Y^{s,r}_n,M})(Y^{s,r}_n),(\Yl(t_n^k)-\Yl(t_n))s\rangle]_{Z=\Yl(t_n)}.
\end{align*}
Thus
\begin{align*}
	\partial_{\mu} &f_m(Z,\Pls)(Y_n^s)=\partial_{\mu} f_m(Z,\Pl)(\Yl(t_n)) \\
	&+\E[ \partial^2_{\mu}f_m(Z,\mathcal L_{h_l}^{\e,Y^{s,r}_n,M})(Y^{s,r}_n),(\Yl(t_n^k)-\Yl(t_n))s]_{Z=\Yl(t_n)}.
\end{align*}
Substituting the last equation into the second summand of the RHS of Equation \eqref{eq1-bar ABE} yields
\begin{align*}
	&f_m(\Yl(t_n),\Plk)-f_m(\Yl(t_n),\Pl) \\
	&=\E[\langle \partial_{\mu} f_m(Z,\Pls)(Y_n^s),\sum_{j=0}^{k-1}f(\Yl(t_n^j),\Plj)h_l \rangle]_{Z=\Yl(t_n)} \\
	&+ \E[\langle \partial_{\mu} f_m(Z,\Pl)(\Yl(t_n)), \e \sqrt{h_l} \sum_{j=0}^{k-1}g(\Yl(t_n^j),\Plj)\D \xi_n^j \rangle]_{Z=\Yl(t_n)} \\
	&+ \E[\langle \partial^2_{\mu}f_m(Z,\mathcal L_{h_l}^{\e,Y^{s,r}_n,M})(Y_n^{s,r})(\Yl(t_n^k)-\Yl(t_n))s, \\
	& \quad \quad \quad \quad \quad \quad \quad \quad \quad \quad \quad \quad \quad \quad \quad \quad \quad \e \sqrt{h_l} \sum_{j=0}^{k-1}g(\Yl(t_n^j),\Plj)\D \xi_n^j \rangle]_{Z=\Yl(t_n)}.
\end{align*}
By independence the second expectation above is zero, therefore the proof is complete.

\hfill $\Box$

{\bf Proof of Theorem \ref{2ndMom2Paths}}
From \eqref{MMC EM scheme_2} and \eqref{MMC EM scheme_3} we have that for $n \leq N^{l-1}-1$ 
\begin{align*}
\Yl(\tnp) &- \Ylm(\tnp)=\Yl(t_n) - \Ylm(t_n) \\
&+h_l \sum_{k=0}^{N-1}\left(f(\Yl(t_n^k),\Plk)-f(\Yl(t_n),\Pl)\right) \\
&+h_l \sum_{k=0}^{N-1}\left(f(\Yl(t_n),\Pl)-f(\Ylm(t_n),\Plm)\right)  \\
&+\e \sqrt{h_l} \sum_{k=0}^{N-1}\left(g(\Yl(t_n^k),\Plk)-g(\Yl(t_n),\Pl)\right)\D \xi^k_n \\
&+\e \sqrt{h_l}\sum_{k=0}^{N-1}\left(g(\Yl(t_n),\Pl)-g(\Ylm(t_n),\Plm)\right)\D \xi^k_n  \\
&=: \Yl(t_n) - \Ylm(t_n) + R_N.
\end{align*}
By using the linearity property of the inner product, we obtain
\begin{align*}
	|\Yl(\tnp) &- \Ylm(\tnp)|^2=\langle \Yl(t_n) - \Ylm(t_n)+R_N, \Yl(t_n) - \Ylm(t_n)+R_N \rangle \\
	&=|\Yl(t_n) - \Ylm(t_n)|^2+|R_N|^2+2\langle \Yl(t_n) - \Ylm(t_n), R_N \rangle.  
\end{align*}
Applying the elementary inequality $|a+b+c+d|^2 \leq 4|a|^2+4|b|^2+4|c|^2+4|d|^2$ to the term $|R_N|^2$ above, we derive that
\begin{align*}
| \Yl&(\tnp) - \Ylm(\tnp) |^2 \leq |\Yl(t_n) - \Ylm(t_n) |^2	\\
&+ 4 h^2_l \left| \sum_{k=0}^{N-1}\left(f(\Yl(t_n^k),\Plk)-f(\Yl(t_n),\Pl)\right) \right|^2 \\
&+ 4 h^2_l \left| \sum_{k=0}^{N-1}\left(f(\Yl(t_n),\Pl)-f(\Ylm(t_n),\Plm)\right) \right|^2 \\
&+4 \e^2 \left|  \sum_{k=0}^{N-1}\left(g(\Yl(t_n^k),\Plk)-g(\Yl(t_n),\Pl)\right)\sqrt{h_l}\D \xi^k_n \right|^2 \\
&+4 \e^2 \left| \sum_{k=0}^{N-1}\left(g(\Yl(t_n),\Pl)-g(\Ylm(t_n),\Plm)\right)\sqrt{h_l}\D \xi^k_n \right|^2  \\
&+ 2 h_l \sum_{k=0}^{N-1}\langle \Yl(t_n) - \Ylm(t_n),f(\Yl(t_n^k),\Plk)-f(\Yl(t_n),\Pl) \rangle \\
&+ 2 h_l \sum_{k=0}^{N-1}\langle \Yl(t_n) - \Ylm(t_n),f(\Yl(t_n),\Pl)-f(\Ylm(t_n),\Plm) \rangle \\
&+ 2 \e \sqrt{h_l} \sum_{k=0}^{N-1}\langle \Yl(t_n) - \Ylm(t_n),\big(g(\Yl(t_n^k),\Plk)-g(\Yl(t_n),\Pl)\big) \D \xi^k_n \rangle \\
&+ 2 \e \sqrt{h_l} \sum_{k=0}^{N-1}\langle \Yl(t_n) - \Ylm(t_n),\big(g(\Yl(t_n),\Pl)-g(\Ylm(t_n),\Plm)\big) \D \xi^k_n \rangle. 
\end{align*}
Now, we take expectations on both sides of the previous inequality. Since $\Delta \xi^k_n$ is independent of $\Yl(t^k_n)$ and $\Ylm(t_n)$, the expectation of the last two summands in the equation above is zero. Thus,
\begin{align}\label{eq1 - theo 2nd moment}
	\E [| \Yl&(\tnp) - \Ylm(\tnp) |^2 ] \leq \E [ |\Yl(t_n) - \Ylm(t_n) |^2] \\
	&+ 4 N h^2_l  \sum_{k=0}^{N-1} \E \left| f(\Yl(t_n^k),\Plk)-f(\Yl(t_n),\Pl) \right|^2 \nonumber \\
	&+ 4 N h^2_l  \sum_{k=0}^{N-1} \E \left| f(\Yl(t_n),\Pl)-f(\Ylm(t_n),\Plm) \right|^2 \nonumber \\
	&+4  \e^2 \E\left[ \left|\sum_{k=0}^{N-1}  \left(g(\Yl(t_n^k),\Plk)-g(\Yl(t_n),\Pl)\right)\sqrt{h_l}\D \xi^k_n \right|^2 \right] \nonumber \\
	&+4 \e^2 \E\left[ \left| \sum_{k=0}^{N-1}    \left(g(\Yl(t_n),\Pl)-g(\Ylm(t_n),\Plm)\right)\sqrt{h_l}\D \xi^k_n \right|^2\right]  \nonumber \\
	&+ 2 h_l \sum_{k=0}^{N-1} \E [\langle \Yl(t_n) - \Ylm(t_n),f(\Yl(t_n^k),\Plk)-f(\Yl(t_n),\Pl) \rangle] \nonumber\\
	&+ 2 h_l \sum_{k=0}^{N-1} \E [\langle \Yl(t_n) - \Ylm(t_n),f(\Yl(t_n),\Pl)-f(\Ylm(t_n),\Plm) \rangle]. \nonumber \\
	&=: \E [ |\Yl(t_n) - \Ylm(t_n) |^2]  + I_1+I_2+I_3+I_4+I_5+I_6.
\end{align}
By Assumption \ref{a1} and Lemma \ref{lemma kn}, one can see that
\begin{align*}
	I_1 &\leq 4 K N  h_l^2 \sum_{k=0}^{N-1} ( \E|\Yl(t_n^k)-\Yl(t_n)|^2 + \mathbb W_2^2(\Plk, \Pl)) \\
	&\leq 8 K N h_l^2 \sum_{k=0}^{N-1}  \E|\Yl(t_n^k)-\Yl(t_n)|^2 \leq 8 K N^2 h^2_l(C N^2 h^2_l+C N\e^2h_l).
\end{align*}
Also, by Assumption \ref{a1}
\begin{align*}
	I_2 &\leq 4 K N  h_l^2 \sum_{k=0}^{N-1}  (\E|\Yl(t_n)-\Ylm(t_n)|^2 + \mathbb W_2^2(\Pl, \Plm) )\\
	&\leq 8 K N^2 h_l^2 \E [ |\Yl(t_n) - \Ylm(t_n) |^2].
\end{align*}
By the BDG inequality, Assumption \ref{a1} and Lemma \ref{lemma kn}, we obtain 
\begin{align*}
I_3 &\leq C \e^2  \sum_{k=0}^{N-1}\E[| g(\Yl(t_n^k),\Plk)-g(\Yl(t_n),\Pl)|^2] h_l\\
&=C  h_l \e^2 \sum_{k=0}^{N-1}(\E[|\Yl(t^k_n)-\Yl(t_n)|^2] + \mathbb W_2^2(\Plk, \Pl)) \leq C N^3 h^3_l \e^2 + C  N^2 h_l^2 \e^4.                   
\end{align*}
Similarly to $I_3$,
$$I_4 \leq CN h_l \e^2 \E[|\Yl(t_n) - \Ylm(t_n)|^2].$$
An application of the Cauchy-Schwarz inequality and Assumption \ref{a1} gives
\begin{align*}
	I_5 &=  2 h_l \sum_{k=0}^{N-1} \E [\langle \Yl(t_n) - \Ylm(t_n),f(\Yl(t_n^k),\Plk)-f(\Yl(t_n),\Pl) \rangle]  \\
	&=  2 h_l \sum_{k=0}^{N-1} \E [\langle \Yl(t_n) - \Ylm(t_n),f(\Yl(t_n^k),\Plk)-f(\Yl(t_n),\Plk) \rangle]  \\
	&+ 2 h_l \sum_{k=0}^{N-1} \E [\langle \Yl(t_n) - \Ylm(t_n),f(\Yl(t_n),\Plk)-f(\Yl(t_n),\Pl) \rangle] \\
	&=: I_{5A} + I_{5B}.
\end{align*}
Applying Lemma \ref{ABE} we have
\begin{align*}
I_{5A} &\leq  2 h_l \sum_{k=0}^{N-1} \E [\langle \Yl(t_n) - \Ylm(t_n),A_k  \rangle] +2 h_l \sum_{k=0}^{N-1} \E [\langle \Yl(t_n) - \Ylm(t_n),B_k  \rangle] \\
&+ 2 h_l \sum_{k=0}^{N-1} \E [\langle \Yl(t_n) - \Ylm(t_n),E_k  \rangle].
\end{align*}
By independence, the second summand above is zero. 
Also, we note that
\begin{align*}
	\E[|A_k|^2] &= \sum_{m=1}^{\bar d}\E[|A^m_k|]^2  \leq \bar d \bar K \E \left|h_l \sum_{j=0}^{k-1}f(\Yl(t^j_n),\Plj)  \right|^2  \\
	&\leq \bar d \bar K h^2_j N \sum_{j=0}^{k-1} \E\left[\left(\beta \left(1+|\Yl(t^j_n)|^2 + W_2^2(\Plj) \right)\right)^2\right] \\
	&\leq  \bar K h_l^2 N^2 C
\end{align*}
and
\begin{align}\label{2ndE}
	\E[|E_k|^2] &= \sum_{m=1}^{\bar d}\E[(E^m_k)^2] \leq \bar d K \e^2h_l \E\left[|\Yl(t^k_n)-\Yl(t_n)|^2 \left|\e \sqrt{h_l} \sum_{j=0}^{k-1} g(\Yl(t^j_n),\Plj) \D \xi_n^j \right|^2 \right] \nonumber \\
	&\leq \bar d \bar K \e^2 h_l ( \E[|\Yl(t^k_n)-\Yl(t_n)|^4])^{1/2}\left(\E\left[\left|\sum_{j=0}^{k-1} g(\Yl(t^j_n),\Plj) \D \xi_n^j \right|^4\right]\right)^{1/2} \nonumber \\
	&\leq  \bar K \e^2 C N^3 h_l^3+  \bar K \e^4 C N^2 h_l^2,
\end{align}
where Lemma \ref{lemma kn} is used in the last inequality. Therefore, applying the Cauchy-Schwartz inequality first and the elementary inequality $2ab \leq a^2 + b^2$ later yields
\begin{align*}
I_{5A} &\leq 2 h_l \sum_{k=0}^{N-1} \E [| \Yl(t_n) - \Ylm(t_n)||A_k | ] +2 h_l \sum_{k=0}^{N-1} \E [| \Yl(t_n) - \Ylm(t_n)||E_k|] \\
&\leq 2 h_l \sum_{k=0}^{N-1} \E [| \Yl(t_n) - \Ylm(t_n)|^2+ h_l \sum_{k=0}^{N-1} \E [|A_k |^2]+ h_l \sum_{k=0}^{N-1} \E [|E_k |^2] \\
&\leq 2 h_l N \E [| \Yl(t_n) - \Ylm(t_n)|^2 +\bar K h_l^3 N^3 C + \bar K  C N^4 h_l^4 \e^2 + \bar K C N^3 h_l^3 \e^4
\end{align*}
	
Similarly, using Lemma \ref{bar ABE} one can see that 

\begin{align*}
I_{5B} &\leq 2 h_l \sum_{k=0}^{N-1} \E [\langle \Yl(t_n) - \Ylm(t_n),\bar A_k  \rangle] +2 h_l \sum_{k=0}^{N-1} \E [\langle \Yl(t_n) - \Ylm(t_n),\bar E_k  \rangle] 
\end{align*}	
Also, we have $\E[|\bar A_k|^2] \leq h_l^2 N^2 C$ and
\begin{equation}\label{2nd bar E}
\E[|\bar E_k|^2] \leq \bar K \e^2 C N^3 h_l^3+\bar K \e^4 C N^2 h_l^2.
\end{equation}
Thus,
\begin{align*}
	I_{5B} &\leq 2 h_l \sum_{k=0}^{N-1} \E [| \Yl(t_n) - \Ylm(t_n)|\E[|\bar A_k |] ] +2 h_l \sum_{k=0}^{N-1} \E [| \Yl(t_n) - \Ylm(t_n)| \E[|\bar E_k|]] \\
	&\leq 2 h_l \sum_{k=0}^{N-1} \E [| \Yl(t_n) - \Ylm(t_n)|^2+ h_l \sum_{k=0}^{N-1} \E [|\bar A_k |^2]+ h_l \sum_{k=0}^{N-1} \E [|\bar E_k |^2] \\
	&\leq 2 h_l N \E [| \Yl(t_n) - \Ylm(t_n)|^2 +\bar K h_l^3 N^3 C + \bar K  C N^4 h_l^4 \e^2 + \bar K C N^3 h_l^3 \e^4
\end{align*}
Additionally, we have 
\begin{align*}
	I_6 &\leq h_l N \E | \Yl(t_n) - \Ylm(t_n)|^2 +h_l N \E | \Yl(t_n) - \Ylm(t_n)|^2 + h_l N \mathbb W^2_2(\Plk,\Plkm) \\
	&\leq 3 h_l N \E | \Yl(t_n) - \Ylm(t_n)|^2.
\end{align*} 
Substituting the bounds for the terms $I_1$ to $I_6$ into Equation \eqref{eq1 - theo 2nd moment} yields that for $n \leq N^{l-1}-1$ 
\begin{align*}
	\E [| \Yl(\tnp) - \Ylm(\tnp) |^2 ] &\leq \E [ |\Yl(t_n) - \Ylm(t_n) |^2] + \hat C\E [ |\Yl(t_n) - \Ylm(t_n) |^2] \\
	&+   C N^3 h_l^3 + C N^2 h_l^2 \e^4, 
\end{align*}
which implies that that for $n \leq N^{l-1}-1$
\begin{align*}
	\E [| \Yl(\tnp) - \Ylm(\tnp) |^2 ] &\leq   \hat C \sum_{k=1}^n \E [ |\Yl(t_k) - \Ylm(t_k) |^2] \\
	&+   C N^2 h_l^2 + C N h_l \e^4.
\end{align*}

An application of the discrete Gronwall inequality yields the result.
	
\hfill $\Box$

\subsection{Estimates of Variance}

In this section we provide an estimate for the variance of two coupled paths which is the main result of the paper and will be presented in  Theorem \ref{th1}. We will need the following lemma taken from \cite{ahs15}. Proof of this theorem can be found in \cite{andersonII}.

\begin{lemma}\label{vl2}
	Suppose that $A^{\e,h}$ and $B^{\e,h}$ are families of random variables determined by scaling parameters $\e$ and $h$. Further, suppose that there are $C_1>0,C_2>0$ and $C_3 >0$ such that for all $\e \in (0,1)$ the following three conditions hold:
	\begin{enumerate}[label=(\roman*)]
	\item ${\rm Var}(A^{\e,h}) \leq C_1\e^2 \ \text{uniformly in $h$}$,
	\item $|A^{\e,h}| \leq C_2 \ \text{uniformly in $h$},$
	\item $|\E[B^{\e,h}]| \leq C_3 h.$
	\end{enumerate}
Then
$${\rm Var}(A^{\e,h} B^{\e,h}) \leq 3C^2_3C_1 h^2 \e^2 + 15 C_2^2 {\rm Var}(B^{\e,h}).$$
\end{lemma}
The following two lemmas that will be needed to prove Theorem \ref{th1}.
\begin{lemma}\label{vl1}
	 Assume that $\gamma:\mathbb R^d \rightarrow \mathbb R$ satisfies the Lipschitz condition, i.e. for all $x,y \in \RR^d$ there exists a positive constant $L$, such that $|\gamma(x)-\gamma(y)|^2 \leq L|x-y|^2.$ Then for $s \in [0,1]$ one has
	$$\max_{\substack{0 \leq n \leq N^{l-1} \\ 1 \leq k \leq N}} {\rm Var}(\gamma(sY^{\e,i,M}_{h_{l_2}}(t^k_n)+(1-s)(Y^{\e,i,M}_{h_{l_1}}(t_n)))) \leq C \e^2.$$
\end{lemma}	
{\bf Proof.}
Let $z_{h_{l_1}}$ and $z_{h_{l_2}}$ be defined by \eqref{z_h}. Using the fact that for a random variable $X$ and a constant $a,$  ${\rm Var}(X+a)={\rm Var}(X)$ and the fact that $\gamma$ is Lipschitz, we have that
\begin{align*}
	&\max_{\substack{0 \leq n \leq N^{l-1} \\ 1 \leq k \leq N}} {\rm Var}(\gamma(sY^{\e,i,M}_{h_{l_2}}(t^k_n)+(1-s)(Y^{\e,i,M}_{h_{l_1}}(t_n)))) \\
	&= \max_{\substack{0 \leq n \leq N^{l-1} \\ 1 \leq k \leq N}} {\rm Var}(\gamma(sY^{\e,i,M}_{h_{l_2}}(t^k_n)+(1-s)(Y^{\e,i,M}_{h_{l_1}}(t_n)))-\gamma(sz_{h_{l_2}}(t^k_n)+(1-s)(z_{h_{l_1}}(t_n)))) \\
	&\leq \max_{\substack{0 \leq n \leq N^{l-1} \\ 1 \leq k \leq N}}\E[|(\gamma(sY^{\e,i,M}_{h_{l_2}}(t^k_n)+(1-s)(Y^{\e,i,M}_{h_{l_1}}(t_n)))-\gamma(sz_{h_{l_2}}(t^k_n)+(1-s)(z_{h_{l_1}}(t_n)))|^2]\\
	&= \max_{\substack{0 \leq n \leq N^{l-1} \\ 1 \leq k \leq N}}L\E[|sY^{\e,i,M}_{h_{l_2}}(t^k_n)+(1-s)(Y^{\e,i,M}_{h_{l_1}}(t_n))-sz_{h_{l_2}}(t^k_n)-(1-s)(z_{h_{l_1}}(t_n))|^2] \\
	&\leq  \max_{\substack{0 \leq n \leq N^{l-1} \\ 1 \leq k \leq N}}sL\E[|(Y^{\e,i,M}_{h_{l_2}}(t^k_n)-z_{h_{l_2}}(t^k_n)|^2]+(1-s)L\E[|(Y^{\e,i,M}_{h_{l_1}}(t_n)-(z_{h_{l_1}}(t_n)|^2].
\end{align*}
The required assertion follows by Lemma \ref{z}.	
\hfill $\Box$

\begin{lemma}\label{cond3}
Let Assumption \ref{a1} hold. Then there exists a positive constant $C$ such that
	$$\max_{\substack{0 \leq n \leq N^{l-1} \\ 1 \leq k \leq N}} |\E[\Yl(t^k_n) - \Yl(t_n)]| \leq CN h_l.$$	
\end{lemma}	
{\bf Proof.} 
From \eqref{MMC EM scheme_1} we  have that
\begin{align*}	
	|\E[\Yl(t^k_n) &- \Yl(t_n)]| \\
	&=\left|  \sum_{j=0}^{k-1}\E[f(\Yl(t_n^j),\Plj)]h_l + \e \sqrt{h_l} \sum_{j=0}^{k-1}\E[g(\Yl(t_n^j),\Plj)\D \xi_n^j]\right|. 
\end{align*}
By independence the second summand of RHS in  above is zero. Thus using Jensen's inequality and Remark \ref{onem} yields
\begin{align*}
	|\E[\Yl(t^k_n) - \Yl(t_n)]|	&\leq \sum_{j=0}^{k-1}\E[|f(\Yl(t_n^j),\Plj)|]h_l \\
	&\leq  h_l \sum_{j=0}^{k-1} \E[\sqrt \beta \big(1+|\Yl(t^j_n)|^2 + W_2^2(\Plj)\big)^{1/2} ] \\
	&\leq  \sqrt \beta h_l \sum_{j=0}^{k-1} \left(1+2\E[|\Yl(t^j_n)|^2] \right)^{1/2}.
\end{align*}
An application of Lemma \ref{0pmoment} and the fact that $k \leq N,$ completes the proof.
\hfill $\Box$

Now, we can formulate the main result of the paper.

\begin{theorem}\label{th1} 
Let Assumption \ref{a1} hold, assume that $\Psi:\mathbb{R}^d\rightarrow\mathbb{R}$ has continuous second order derivative and there exists a constant $C$ such that
\begin{equation*}
\begin{split}
\left|\frac{\partial \Psi}{\partial x_i}\right|\le C~~and ~~\left|\frac{\partial^2 \Psi}{\partial x_i\partial x_j}\right|\le C
\end{split}
\end{equation*}
for any $i,j=1,2,\cdots,a$. Then, we have
\begin{equation*}
\begin{split}
\max_{0\le n<M^{l-1}}{\rm Var}(\Psi(Y_{h_{l}}^{\vv, i, M}(t_{n+1}))-\Psi(Y_{h_{l-1}}^{\vv, i, M}(t_{n+1}))\le  C\vv^2 h_{l-1}^2+C\vv^4h_{l-1}.
\end{split}
\end{equation*}
\end{theorem}

{\bf Proof.} 
From \eqref{MMC EM scheme_2} and \eqref{MMC EM scheme_3} we have that for $n \leq N^{l-1}-1$ 
\begin{align*}
	[\Yl(\tnp) &- \Ylm(\tnp)]_j=[\Yl(t_n) - \Ylm(t_n)]_j \\
	&+h_l \sum_{k=0}^{N-1}\left(f_j(\Yl(t_n^k),\Plk)-f_j(\Yl(t_n),\Pl)\right) \\
	&+h_l \sum_{k=0}^{N-1}\left(f_j(\Yl(t_n),\Pl)-f_j(\Ylm(t_n),\Plm)\right)  \\
	&+\e \sqrt{h_l} \sum_{k=0}^{N-1}\left(g_j(\Yl(t_n^k),\Plk)-g_j(\Yl(t_n),\Pl)\right)\D \xi^k_n \\
	&+\e \sqrt{h_l}\sum_{k=0}^{N-1}\left(g_j(\Yl(t_n),\Pl)-g_j(\Ylm(t_n),\Plm)\right)\D \xi^k_n,
\end{align*}
where $f_j$ is the $j$th component of $f$ and $g_j$ is the $j$th row of $g.$ Taking variances on both sides of the previous inequality and using simple properties of variance and covariance functions, we obtain
\begin{align*}
	{\rm Var}&([\Yl(\tnp) - \Ylm(\tnp)]_j)\leq(1+N h_l){\rm Var}([\Yl(t_n) - \Ylm(t_n)]_j) \\
	&+4 h^2_l N \sum_{k=0}^{N-1}{\rm Var}\left(f_j(\Yl(t_n^k),\Plk)-f_j(\Yl(t_n),\Pl)\right) \\
	&+(4N h_l + 1)N h_l {\rm Var}\left(f_j(\Yl(t_n),\Pl)-f_j(\Ylm(t_n),\Plm)\right)  \\
	&+4\e^2 h_l \sum_{k=0}^{N-1}{\rm Var}\left(g_j(\Yl(t_n^k),\Plk)-g_j(\Yl(t_n),\Pl)\right)\D \xi^k_n \\
	&+4 \e^2 h_l\sum_{k=0}^{N-1}{\rm Var}\left(g_j(\Yl(t_n),\Pl)-g_j(\Ylm(t_n),\Plm)\right)\D \xi^k_n \\
	&+2 {\rm Cov}\left([\Yl(t_n) - \Ylm(t_n)]_j,h_l \sum_{k=0}^{N-1}f_j(\Yl(t_n^k),\Plk)-f_j(\Yl(t_n),\Pl)\right)\\
	&=: I_1 + I_2 + I_3 + I_4 + I_5 + I_6.
\end{align*}
In order to complete the proof of the theorem, we give estimates for $I_i, i=2,...,6,$ which will be shown in the following lemmas.
\begin{lemma} \label{lem2} There exists a positive constant $C$  such that
$$I_2 \leq  C N^3 h^3_l \e^2.$$
\end{lemma}	
{\bf Proof.} Using the fact that for two random variables $X,Y, {\rm Var}(X+Y)\leq 2{\rm Var}(X)+2{\rm Var}(Y)$, we have that 
\begin{align*}
	{\rm Var}(f_j(&\Yl(t_n^k),\Plk)-f_j(\Yl(t_n),\Pl))\\
	&\leq  2{\rm Var}(f_j(\Yl(t_n^k),\Plk)-f_j(\Yl(t_n),\Plk)) \\
	&+ 2{\rm Var}(f_j(\Yl(t_n),\Plk)-f_j(\Yl(t_n),\Pl))=:I_{2A}+I_{2B}.
\end{align*}
First we estimate $I_{2A}$. By the mean value theorem there exists an $s \in [0,1]$ such that
\begin{align*} 
	f_m(\Yl &(t_n^k),\Plk)-f_j(\Yl(t_n),\Plk)  \\
	&= \langle \nabla f_j(s \Yl(t_n^k) + (1-s)\Yl(t_n)),\Plk), (\Yl(t_n^k) - \Yl(t_n)) \rangle.  \nonumber
\end{align*}
Let $\nabla_q f_j(s \Yl(t_n^k) + (1-s)\Yl(t_n)),\Plk)$ and $[(\Yl(t_n^k) - \Yl(t_n))]_q$ be the $q$ components of $\nabla f_j(s \Yl(t_n^k) + (1-s)\Yl(t_n)),\Plk)$ and $(\Yl(t_n^k) - \Yl(t_n))$ respectively. We want to apply Lemma \ref{vl2} with $A^{\e,h} = \nabla_q f_j(s \Yl(t_n^k) + (1-s)\Yl(t_n)),\Plk)$ and $B^{\e,h}=[(\Yl(t_n^k) - \Yl(t_n))]_q$ so we check that the three conditions are satisfied. By Assumption \ref{a1}, the function $\nabla^2_q f_j$ is bounded, so $\nabla_q f_j$ is Lipschitz on the first argument. Applying Lemma \ref{vl1} with $\gamma=\nabla_q f_j(\cdot,\Plk)$ and $h_{l_1} = h_{l_2} = h_l,$ we obtain
\begin{equation}\label{1stcond} 
{\rm Var}(\nabla_q f_j(s \Yl(t_n^k) + (1-s)\Yl(t_n)),\Plk)) \leq C_1 \e^2,
\end{equation}
so the first condition of Lemma \ref{vl2} is satisfied. Conditions 2 and 3 are satisfied by Assumption \ref{a1} and Lemma \ref{cond3} respectively. Thus by Lemma \ref{vl2} we have that
\begin{align*}
{\rm Var}(\nabla_q f_j(s \Yl(t_n^k) &+ (1-s)\Yl(t_n)),\Plk)[(\Yl(t_n^k) - \Yl(t_n))]_q) \\
 &\leq 3C^2_3C_1 N^2 h_l^2 \e^2 + 15 C_2^2 {\rm Var}([(\Yl(t_n^k) - \Yl(t_n))]_q).
\end{align*}
In order to estimate ${\rm Var}([(\Yl(t_n^k) - \Yl(t_n))]_q)$ we use Equation \eqref{MMC EM scheme_1} to obtain
\begin{align*}	
	{\rm Var}&([(\Yl(t_n^k) - \Yl(t_n))]_q) \\
	 &\leq 2{\rm Var} (\sum_{j=0}^{k-1}f_q(\Yl(t_n^j),\Plj)h_l) + 	2{\rm Var}(\e \sqrt{h_l} \sum_{j=0}^{k-1}g_q(\Yl(t_n^j),\Plj)\D \xi_n^j). 
\end{align*}
By Asumption \ref{a1} and Lemma \ref{z} we have that
\begin{align*}	
   {\rm Var} (\sum_{j=0}^{k-1}&f_q(\Yl(t_n^j),\Plj)h_l) = {\rm Var} (h_l\sum_{j=0}^{k-1}f_q(\Yl(t_n^j),\Plj)-f_q(z_h(t_n^j),\delta_{z_h(t_n^j)})) \\
   &\leq h_l^2 \E [|(\sum_{j=0}^{k-1}f_q(\Yl(t_n^j),\Plj)-f_q(z_h(t_n^j),\delta_{z_h(t_n^j)}))|^2] \leq C N^2 h_l^2 \e^2.
\end{align*}
From \eqref{eq3-lemma1} we have that
$${\rm Var}(\e \sqrt{h_l} \sum_{j=0}^{k-1}g_q(\Yl(t_n^j),\Plj)\D \xi_n^j) \leq C N h_l \e^2.$$
Thus
$${\rm Var}([(\Yl(t_n^k) - \Yl(t_n))]_q) \leq  C N^2 h_l^2 \e^2 +  C N h_l \e^2.$$
Using the formula ${\rm Var} (\sum_{i=1}^d X_i )\leq d \sum_{i=1}^d {\rm Var}(X_i)$ with $i=q,X_i = [\Yl(t_n^k) - \Yl(t_n)]_q$ yields
$${\rm Var}([(\Yl(t_n^k) - \Yl(t_n))]_q) \leq  d^2C N^2 h_l^2 \e^2 + d^2  C N h_l \e^2 \leq C N h_l \e^2.$$
Thus,
\begin{align*}
	I_{2A} \leq  C N h_l \e^2 .
\end{align*}
Next, we estimate $I_{2B}$. By Equation \eqref{mvt} there exists a random variable $s:\Omega \rightarrow [0,1]$ such that
\begin{align*}
	f_j(\Yl(t_n),\Plk)-&f_j(\Yl(t_n),\Pl) \\
	&=\E[\langle \partial_{\mu} f_j(Z,\Pls)(Y_n^s),(\Yl(t^k_n)-\Yl(t_n)) \rangle]_{Z=\Yl(t_n)}.
\end{align*}
where $Y_n^s := s \Yl(t_n^k) + (1-s)\Yl(t_n).$  Let $\partial_{\mu,q}f_j(Z,\Pls)(Y_n^s)$ and $[\Yl(t^k_n)-\Yl(t_n)]_q$ be the $q$-components of $\partial_{\mu}f_j(Z,\Pls)(Y_n^s)$ and $\Yl(t^k_n)-\Yl(t_n)$ respectively. Then
\begin{align*}
	&{\rm Var}(\E[\partial_{\mu,q}f_j(Z,\Pls)(Y_n^s)[\Yl(t^k_n)-\Yl(t_n)]_q]_{Z=\Yl(t_n)}) \\
	&={\rm Var}(\E[\partial_{\mu,q}f_j(Z,\Pls)(Y_n^s)[\Yl(t^k_n)-\Yl(t_n)]_q]_{Z=\Yl(t_n)}\\
	&-\E[\partial_{\mu,q}f_j(z_{h_l}(t_n),\delta_{z_{h_l}(t_n)})(z_{h_l}(t_n))[\Yl(t^k_n)-\Yl(t_n)]_q] ) \\
	&=	{\rm Var}(\E[(\partial_{\mu,q}f_j(Z,\Pls)(Y_n^s)-\partial_{\mu,q}f_j(z_{h_l}(t_n),\delta_{z_{h_l}(t_n)})(z_{h_l}(t_n)))\\
	&\times [\Yl(t^k_n)-\Yl(t_n)]_q]]_{Z=\Yl(t_n)} )\\
	&\leq \E[(\E[(\partial_{\mu,q}f_j(Z,\Pls)(Y_n^s)-\partial_{\mu,q}f_j(z_{h_l}(t_n),\delta_{z_{h_l}(t_n)})(z_{h_l}(t_n)))\\
	&\times [\Yl(t^k_n)-\Yl(t_n)]_q]]_{Z=\Yl(t_n)} )^2] \\
	&\leq \E[\E[|\partial_{\mu,q}f_j(Z,\Pls)(Y_n^s)-\partial_{\mu,q}f_j(z_{h_l}(t_n),\delta_{z_{h_l}(t_n)})(z_{h_l}(t_n))|^2]_{Z=\Yl(t_n)}\\
	&\times\E[|[\Yl(t^k_n)-\Yl(t_n)]_q|^2]],\\
\end{align*}
where we have use the Cauchy-Schwarz inequality in the penultimate step.
By condition \eqref{a2a} and Lemma \ref{z}
$$\E[\E[|\partial_{\mu,q}f_j(Z,\Pls)(Y_n^s)-\partial_{\mu,q}f_j(z_{h_l}(t_n),\delta_{z_{h_l}(t_n)})(z_{h_l}(t_n))|^2]_{Z=\Yl(t_n)} \leq C\e^2$$
and by Lemma \ref{lemma kn}
$$\E[|\Yl(t^k_n) - \Yl(t_n)|^2] \leq C N^2 h_l^2 +C N   h_l \e^2.$$
Therefore
$$I_{2B} \leq C N^2 h^2_l \e^2 + C N h_l \e^4,$$
and the proof is complete.
\hfill $\Box$

\begin{lemma}
	There exists positive constants $C$ and $\bar C$ such that
	$$I_3 \leq C N h_l \sum_{q=1}^d {\rm Var}([(\Yl(t_n) - \Ylm(t_n))]_q +C N^3 h_l^3 \e^2.$$
\end{lemma}	
{\bf Proof.} Note that 
\begin{align*}
	{\rm Var}(f_j&(\Yl(t_n),\Pl)-f_j(\Ylm(t_n),\Plm))\\
	&\leq  2{\rm Var}(f_j(\Yl(t_n),\Pl)-f_j(\Ylm(t_n),\Pl)) \\
	&+ 2{\rm Var}(f_j(\Ylm(t_n),\Pl)-f_j(\Ylm(t_n),\Plm))=:I_{3A} + I_{3B}.
\end{align*}
First, we estimate $I_{3B}.$ By the mean value theorem there exists an $s \in [0,1]$ such that
\begin{align*} 
	f_j(\Yl &(t_n),\Pl)-f_j(\Ylm(t_n),\Pl)  \\
	&= \langle \nabla f_j(s \Yl(t_n) + (1-s)\Ylm(t_n)),\Pl), (\Yl(t_n) - \Ylm(t_n)) \rangle.  \nonumber
\end{align*}
Let $\nabla_q f_j(s \Yl(t_n) + (1-s)\Ylm(t_n)),\Pl)$ and $[(\Yl(t_n) - \Ylm(t_n))]_q$ be the $q$ components of $\nabla f_j(s \Yl(t_n) + (1-s)\Ylm(t_n)),\Pl)$ and $(\Yl(t_n) - \Ylm(t_n))$ respectively. We want to apply Lemma \ref{vl2} with $A^{\e,h} = \nabla_q f_j(s \Yl(t_n) + (1-s)\Ylm(t_n)),\Pl)$ and $B^{\e,h}=[(\Yl(t_n) - \Ylm(t_n))]_q$ so we check that the three conditions are satisfied. Applying Lemma \ref{vl1} with $\gamma=\nabla_q f_j, k=0,h_{l_1}=h_{l-1} $ and $h_{l_2} = h_l,$ we obtain 
$${\rm Var}(\nabla_q f_j(s \Yl(t_n) + (1-s)\Ylm(t_n)),\Pl)) \leq C_1 \e^2,$$ so the first condition of Lemma \ref{vl2} is satisfied. Conditions 2 and 3 are satisfied by Assumption \ref{a1} and Lemma \ref{cond3} respectively. Thus by Lemma \ref{vl2} we have that
\begin{align*}
	{\rm Var}(\nabla_q f_j(s \Yl(t_n) &+ (1-s)\Ylm(t_n)),\Pl)[(\Yl(t_n) - \Ylm(t_n))]_q) \\
	&\leq 3C^2_3C_1 N^2 h_l^2 \e^2 + 15 C_2^2 {\rm Var}([(\Yl(t_n) - \Ylm(t_n))]_q).
\end{align*}
Using the formula ${\rm Var} (\sum_{i=1}^d X_i )\leq d \sum_{i=1}^d {\rm Var}(X_i)$ with $i=q,X_i = [\Yl(t_n) - \Ylm(t_n)]_q$ yields
$${\rm Var}((\Yl(t_n) - \Ylm(t_n))) \leq C \sum_{q=1}^d {\rm Var}([(\Yl(t_n) - \Ylm(t_n))]_q +C N^2 h_l^2 \e^2.$$
Therefore,
$$I_{3A} \leq CN^2h_l^2\e^2+ C \sum_{q=1}^d {\rm Var}([(\Yl(t_n) - \Ylm(t_n)]_q).$$
Next we estimate $I_{3B}.$ By Equation \eqref{mvt} there exists a random variable $s:\Omega \rightarrow [0,1]$ such that
\begin{align*}
	f_j(\Ylm(t_n),\Pl)-&f_j(\Ylm(t_n),\Plm)  \\
	&=\E[\langle \partial_{\mu} f_j(Z,\Pls)(Y_n^s),(\Yl(t_n)-\Ylm(t_n)) \rangle]_{Z=\Ylm(t_n)}.
\end{align*}
where $Y_n^s := s \Yl(t_n) + (1-s)\Ylm(t_n).$ Let $\partial_{\mu,q}f_j(Z,\Pls)(Y_n^s)$ and $[\Yl(t^n)-\Ylm(t_n)]_q$ be the $q$-components of $\partial_{\mu}f_j(Z,\Pls)(Y_n^s)$ and $\Yl(t_n)-\Ylm(t_n)$ respectively. Then
\begin{align*}
	&{\rm Var}(\E[\partial_{\mu,q}f_j(Z,\Pls)(Y_n^s)[\Yl(t_n)-\Ylm(t_n)]_q]_{Z=\Ylm(t_n)}) \\
	&={\rm Var}(\E[\partial_{\mu,q}f_j(Z,\Pls)(Y_n^s)[\Yl(t^k_n)-\Yl(t_n)]_q]_{Z=\Ylm(t_n)}\\
	&-\E[\partial_{\mu,q}f_j(z_{h_{l-1}}(t_n),\delta_{z_{h_{l-1}}(t_n)})(z_{h_{l-1}}(t_n))[\Yl(t_n)-\Ylm(t_n)]_q] ) \\
	&=	{\rm Var}(\E[(\partial_{\mu,q}f_j(Z,\Pls)(Y_n^s)-\partial_{\mu,q}f_j(z_{h_{l-1}}(t_n),\delta_{z_{h_{l-1}}(t_n)})(z_{h_{l-1}}(t_n)))\\
	&\times[\Yl(t_n)-\Ylm(t_n)]_q]]_{Z=\Ylm(t_n)} )\\
	&\leq \E[(\E[(\partial_{\mu,q}f_j(Z,\Pls)(Y_n^s)-\partial_{\mu,q}f_j(z_{h_{l-1}}(t_n),\delta_{z_{h_{l-1}}(t_n)})(z_{h_{l-1}}(t_n)))\\
	&\times[\Yl(t_n)-\Ylm(t_n)]_q]]_{Z=\Ylm(t_n)} )^2] \\
	&\leq \E[\E[|\partial_{\mu,q}f_j(Z,\Pls)(Y_n^s)-\partial_{\mu,q}f_j(z_{h_{l-1}}(t_n),\delta_{z_{h_{l-1}}(t_n)})(z_{h_{l-1}}(t_n))|^2]_{Z=\Ylm(t_n)}\\
	&\times \E[|[\Yl(t_n)-\Ylm(t_n)]_q|^2]],
\end{align*}
where we have use the Cauchy-Schwarz inequality in the penultimate step.
By condition \eqref{a2a} and Lemma \ref{z}
$$\E[\E[|\partial_{\mu,q}f_j(Z,\Pls)(Y_n^s)-\partial_{\mu,q}f_j(z_{h_{l-1}}(t_n),\delta_{z_{h_{l-1}}(t_n)})(z_{h_{l-1}}(t_n))|^2]_{Z=\Ylm(t_n)} \leq C\e^2$$
and by Theorem \ref{2ndMom2Paths}
$$\E[|\Yl(t_n) - \Ylm(t_n)|^2] \leq C N^2 h_l^2+ C\vv^4 N h_l.$$
Therefore,
$$I_{3B} \leq C N^2 h_l^2 \e^2+ C\vv^6 N h_l,$$
and the proof is complete.

\hfill $\Box$

\begin{lemma}
	There exists a positive constant $C$ such that
	$$I_4 \leq C \e^2 h^3_{l-1} + C \e^4 h^2_{l-1}.$$
\end{lemma}
{\bf Proof.} By Lemma \ref{lemma kn} and Assumption \ref{a1} one can see that
	\begin{align*}
		I_4 &\leq 4\e^2 h_l \sum_{k=0}^{N-1} \E [|g_j(\Yl(t_n^k),\Plk)-g_j(\Yl(t_n),\Pl)|^2] \\
		&\leq 8\e^2 h_l N K (C h_{l-1}^2+  C\e^2 h_{l-1}) = C \e^2 h^3_{l-1} + C \e^4 h^2_{l-1}.
	\end{align*}
\hfill $\Box$

\begin{lemma}
	There exists a positive constant $C$ such that
	$$I_5 \leq C\e^2 h_{l-1}^3+ C \e^6 h^2_{l-1}.$$
\end{lemma}
{\bf Proof.} By Assumption \ref{a1} and Theorem \ref{2ndMom2Paths} we have that
\begin{align*}
	I_5 &\leq 4\e^2 h_l  \sum_{k=0}^{N-1} \E [|g_j(\Yl(t_n),\Pl)-g_j(\Ylm(t_n),\Plm)|^2] \\
	&\leq 4\e^2 h_l N K (C h_{l-1}^2+ C\vv^4 h_{l-1})=C\e^2 h_{l-1}^3+C \e^6 h^2_{l-1}.
\end{align*}
 \hfill $\Box$

\begin{lemma}\label{lem6}
	There exists a positive constant $C$ such that
$$I_6 \leq 2Nh_l {\rm Var}\big([\Yl(t_n) - \Ylm(t_n)]_j\big)+ C N^3 h_l^3 \e^2.$$
\end{lemma}
{\bf Proof.}
Since the covariance is a linear function, by subtracting and adding $f(\Yl(t_n),\Plk)$ to $f_j(\Yl(t_n^k),\Plk)-f_j(\Yl(t_n),\Pl)$ we have that
\begin{align*}
	I_6 &= 2 {\rm Cov}\left([\Yl(t_n) - \Ylm(t_n)]_j,h_l \sum_{k=0}^{N-1}[f_j(\Yl(t_n^k),\Plk)-f_j(\Yl(t_n),\Plk)]\right) \\
	&+ 2 {\rm Cov}\left([\Yl(t_n) - \Ylm(t_n)]_j,h_l \sum_{k=0}^{N-1}[f_j(\Yl(t_n),\Plk)-f_j(\Yl(t_n),\Pl)]\right) \\
	&=: I_{6A} + I_{6B}. 
\end{align*}
By Lemma \ref{ABE}, we obtain
\begin{align*}
	I_{6A} &= 2 {\rm Cov}\left([\Yl(t_n) - \Ylm(t_n)]_j,h_l  \sum_{k=0}^{N-1}(A^j_k + B^j_k+ E^j_k)\right) 	
\end{align*}
Using the bilinearity property of the covariance function we have
\begin{align*}
	I_{6A} &=2h_l \sum_{k=0}^{N-1}{\rm Cov}\Big([\Yl(t_n) - \Ylm(t_n)]_j, A^j_k\Big)+2h_l\sum_{k=0}^{N-1} {\rm Cov}\Big([\Yl(t_n) - \Ylm(t_n)]_j, B^j_k\Big) \\
	&+ 2h_l\sum_{k=0}^{N-1} {\rm Cov}\Big([\Yl(t_n) - \Ylm(t_n)]_j, E^j_k\Big).
\end{align*}
Using the definition of covariance and since the increments $\xi_n^j$ in $B_k^j$ are independent, we find that
\begin{align*}
{\rm Cov}\Big([\Yl(t_n) &- \Ylm(t_n)]_j, B^j_k\Big)\\ &=\E[[\Yl(t_n) - \Ylm(t_n)]_j B^j_k]-\E[[\Yl(t_n) - \Ylm(t_n)]_j]\E[B^j_k]=0.
\end{align*}
Then using the fact that ${\rm Cov}(X,Y) \leq \ff 1 2 {\rm Var}(X) + \ff 1 2 {\rm Var}(Y)$, yields
\begin{equation}\label{I_6var}
I_{6A} \leq 2Nh_l{\rm Var}\big([\Yl(t_n) - \Ylm(t_n)]_j\big)+h_l\sum_{k=0}^{N-1}{\rm Var}(A^j_k)+h_l\sum_{k=0}^{N-1}{\rm Var}(E^j_k).
\end{equation}
Recall from Lemma \ref{ABE} that
$$A^j_k =  \langle \nabla f_j(s \Yl(t_n^k) + (1-s)\Yl(t_n)),\Plk), h_l \sum_{r=0}^{k-1} f(\Yl(t^r_n),\Plr) \rangle.$$
In order to estimate ${\rm Var}(A^j_k)$ we use Lemma \ref{vl2} with $A^{\e,h} = \nabla_q f_j(s \Yl(t_n^k) + (1-s)\Yl(t_n)),\Plk)$ and $B^{\e,h}=[h_l \sum_{r=0}^{k-1} f(\Yl(t^r_n),\Plr)]_q$ so we check that the three conditions are satisfied. The first and second conditions are satisfied by \eqref{1stcond} and Assumption \ref{a1} respectively. By Lemma \ref{0pmoment} and Assumption \ref{a1} we have that
\begin{align*}
	|\E[[h_l \sum_{r=0}^{k-1} f(\Yl(t^r_n),\Plr)]_q]|\leq CNh_l,
\end{align*}
so the third condition is also satisfied. Thus Lemma \ref{vl2} implies that 
$${\rm Var}(A^j_k) \leq C N^2 h^2 \e^2 + C {\rm Var}([h_l \sum_{r=0}^{k-1} f(\Yl(t^r_n),\Plr)]_q).$$
Lemma \ref{z} yields
\begin{align*}
	{\rm Var}([h_l \sum_{r=0}^{k-1} &f(\Yl(t^r_n),\Plr)]_q)={\rm Var}([h_l \sum_{r=0}^{k-1} \{f(\Yl(t^r_n),\Plr)-f(z_{h_l}(t_n^r),\delta_{z_{h_l}(t_n^r)})\}]_q)  \\
	&\leq \E[|([h_l \sum_{r=0}^{k-1} \{f(\Yl(t^r_n),\Plr)-f(z_{h_l}(t_n^r),\delta_{z_{h_l}(t_n^r)})\}]_q)|^2]\leq C N^2 h_l^2 \e^2.
\end{align*}
Therefore
\begin{equation}\label{varA}	
{\rm Var}(A^j_k) \leq C N^2 h^2 \e^2 +C N^2 h_l^2 \e^2 .
\end{equation}
From \eqref{2ndE} we have
\begin{equation}\label{varE}
{\rm Var}(E^j_k)\leq \E[|E^j_k|^2] \leq C N^3 h_l^3 \e^2 +  C N^2 h_l^2 \e^4.
\end{equation}
Substituting \eqref{varA} and \eqref{varE} into \eqref{I_6var} we obtain 
$$I_{6A} \leq 2Nh_l {\rm Var}\big([\Yl(t_n) - \Ylm(t_n)]_j\big)+ C N^3 h_l^3 \e^2.$$
Using Lemma \ref{bar ABE} and simple properties of the covariance function, yields
\begin{align*}
I_{6B} &= 2 {\rm Cov}\left([\Yl(t_n) - \Ylm(t_n)]_j,h_l \sum_{k=0}^{N-1}( \bar A^j_k + \bar E^j_k) \right) \\
&\leq 2h_l \sum_{k=0}^{N-1}{\rm Cov}\Big([\Yl(t_n) - \Ylm(t_n)]_j, \bar A^j_k\Big) + 2h_l\sum_{k=0}^{N-1} {\rm Cov}\Big([\Yl(t_n) - \Ylm(t_n)]_j, \bar E^j_k\Big)\\
&\leq 2Nh_l{\rm Var}\big([\Yl(t_n) - \Ylm(t_n)]_j\big)+h_l\sum_{k=0}^{N-1}{\rm Var}(\bar A^j_k)+h_l\sum_{k=0}^{N-1}{\rm Var}(\bar E^j_k).
\end{align*} 	
Recall from Lemma \ref{bar ABE} that
$$\bar A^j_k =\E[\langle \partial_{\mu} f_j(Z,\Pls)(Y_n^s), h_l \sum_{r=0}^{k-1} f(\Yl(t^r_n),\Plr) \rangle]_{Z=\Yl(t_n)}.$$
Let $\partial_{\mu,q}f_j(Z,\Pls)(Y_n^s)$ and $f_q(\Yl(t^r_n),\Plr)$ be the the $q$-components of $\partial_{\mu} f_j(Z,\Pls)(Y_n^s)$ and $f(\Yl(t^r_n),\Plr)$ respectively. Then
\begin{align*}
{\rm Var}&(\E[\partial_{\mu,q}f_j(Z,\Pls)(Y_n^s)h_l \sum_{r=0}^{k-1} f_q(\Yl(t^r_n),\Plr)]_{Z=\Yl(t_n)}) \\
&={\rm Var}(\E[\partial_{\mu,q}f_j(Z,\Pls)(Y_n^s)h_l \sum_{r=0}^{k-1} f_q(\Yl(t^r_n),\Plr)]_{Z=\Yl(t_n)}\\
&-\E[\partial_{\mu,q}f_j(z_{h_l}(t_n),\delta_{z_{h_l}(t_n)})(z_{h_l}(t_n))h_l \sum_{r=0}^{k-1} f_q(\Yl(t^r_n),\Plr)] ) \\
&=	{\rm Var}(\E[(\partial_{\mu,q}f_j(Z,\Pls)(Y_n^s)-\partial_{\mu,q}f_j(z_{h_l}(t_n),\delta_{z_{h_l}(t_n)})(z_{h_l}(t_n)))\\
&\times h_l \sum_{r=0}^{k-1} f_q(\Yl(t^r_n),\Plr)]]_{Z=\Yl(t_n)} )\\
&\leq \E[(\E[(\partial_{\mu,q}f_j(Z,\Pls)(Y_n^s)-\partial_{\mu,q}f_j(z_{h_l}(t_n),\delta_{z_{h_l}(t_n)})(z_{h_l}(t_n)))\\
&\times h_l \sum_{r=0}^{k-1} f_q(\Yl(t^r_n),\Plr)]_{Z=\Yl(t_n)} )^2] \\
&\leq \E[\E[|\partial_{\mu,q}f_j(Z,\Pls)(Y_n^s)-\partial_{\mu,q}f_j(z_{h_l}(t_n),\delta_{z_{h_l}(t_n)})(z_{h_l}(t_n))|^2]_{Z=\Yl(t_n)}\\
&\times\E[|h_l \sum_{r=0}^{k-1} f_q(\Yl(t^r_n),\Plr)|^2]],
\end{align*}
where we have used the Cauchy-Schwarz inequality in the last step. By condition \eqref{a2a} and Lemma \ref{z}
$$\E[\E[|\partial_{\mu,q}f_j(Z,\Pls)(Y_n^s)-\partial_{\mu,q}f_j(z_{h_l}(t_n),\delta_{z_{h_l}(t_n)})(z_{h_l}(t_n))|^2]_{Z=\Yl(t_n)} \leq C\e^2$$
and by Lemma \ref{0pmoment} and Remark \ref{onem}
$$\E[|h_l \sum_{r=0}^{k-1} f_q(\Yl(t^r_n),\Plr)|^2]] \leq C N^2 h_l^2.$$
Thus,
$${\rm Var}(\bar A^j_k ) \leq C N^2 h_l^2 \e^2.$$
From \eqref{2nd bar E} we have
$${\rm Var}(\bar E^j_k) \leq \E[|\bar E^j_k|^2] \leq \bar K \e^2 C_1 N^3 h_l^3+\bar K \e^4 C N^2 h_l^2.$$
Therefore,
$$I_{6B} \leq 2Nh_l {\rm Var}\big([\Yl(t_n) - \Ylm(t_n)]_j\big)+ C N^3 h_l^3 \e^2$$
and the proof is complete.
\hfill $\Box$

{\bf Continuation of the proof of Theorem \ref{th1}} 
By Lemmas \ref{lem2}-\ref{lem6}, we have
\begin{align*}
	{\rm Var}&([\Yl(\tnp) - \Ylm(\tnp)]_j)\leq{\rm Var}([\Yl(t_n) - \Ylm(t_n)]_j) \\
	&+ C N h_l \sum_{q=1}^d  {\rm Var}([\Yl(t_n) - \Ylm(t_n)]_q)+ C N^3 h_l^3 \e^2 + C N^2 h_l^2 \e^4. 
\end{align*}
Taking the maximum in both sides yields that for  $n \leq N^{l-1}-1$
\begin{align*}
	\max_{1 \leq j \leq d}{\rm Var}&([\Yl(\tnp) - \Ylm(\tnp)]_j)=	\max_{1 \leq j \leq d}{\rm Var}([\Yl(t_n) - \Ylm(t_n)]_j) \\
	&+ C N h_l	\max_{1 \leq j \leq d} {\rm Var}([\Yl(t_n) - \Ylm(t_n)]_j)+ C N^3 h_l^3 \e^2 + C N^2 h_l^2 \e^4. 
\end{align*}
An application of the Grownwall inequality produces
\begin{equation}\label{varY}
\max_{\substack{0 \leq n \leq N^{l-1} \\ 1 \leq j \leq N}}{\rm Var}([\Yl(t_n) - \Ylm(t_n)]_j) \leq C N^2 h_l^2 \e^2 + C N h_l \e^4.
\end{equation}
In order to estimate ${\rm Var}(\Psi(Y_{h_{l}}^{\vv, i, M}(t_n))-\Psi(Y_{h_{l-1}}^{\vv, i, M}(t_n))$ we apply the mean value theorem, so there exists $s \in [0,1]$ such that
$$\Psi(Y_{h_{l}}^{\vv, i, M}(t_n))-\Psi(Y_{h_{l-1}}^{\vv, i, M}(t_n) = \nabla \Psi(s \Yl(t_n) + (1-s)\Ylm(t_n))(\Yl(t_n)-\Ylm(t_n)).$$
We shall apply Lemma \ref{vl2} with $A^{\e,h}=\nabla_q \Psi(s \Yl(t_n) + (1-s)\Ylm(t_n))$ and $B^{\e,h}=[(\Yl(t_n)-\Ylm(t_n))]_q.$
Applying Lemma \ref{vl1} with $\gamma=\nabla_q \Psi, k=0,h_{l_1}=h_{l-1} $ and $h_{l_2} = h_l,$ we obtain 
$${\rm Var}(\nabla_q \Psi(s \Yl(t_n) + (1-s)\Ylm(t_n))) \leq C \e^2,$$ so the first condition of Lemma \ref{vl2} is satisfied. Conditions 2 and 3 are satisfied by Assumption \ref{a1} and Lemma \ref{cond3} respectively. Thus by Lemma \ref{vl2} we have that
\begin{align*}
	{\rm Var}(\nabla_q \Psi(s \Yl(t_n) &+ (1-s)\Ylm(t_n)))[(\Yl(t_n) - \Ylm(t_n))]_q) \\
	&\leq C N^2 h_l^2 \e^2 + C {\rm Var}([(\Yl(t_n) - \Ylm(t_n))]_q).
\end{align*}
Thus
\begin{equation}\label{finVar}
{\rm Var}(\Psi(Y_{h_{l}}^{\vv, i, M}(t_n))-\Psi(Y_{h_{l-1}}^{\vv, i, M}(t_n)) \leq C N^2 h_l^2 \e^2 + C {\rm Var}((\Yl(t_n) - \Ylm(t_n)).
\end{equation}
Sustituting \eqref{varY} into \eqref{finVar} we obtain the desire result.

\hfill $\Box$

\section{Summary}
Regarding the problem of computing $\E[\Phi(X_T)]$ where $X_T$ is the solution at time $T$ to an MV-SDEs with small noise, we studied the problem of comparing the computational cost of using the standard Monte Carlo method with a customatized discretization method versus using the multilevel Monte Carlo method combined with the Euler-Maruyama scheme. To this end, the crucial  part  is to estimate the variance of two coupled paths.  We found that this variance is  $\mathcal O (\vv^2 h_{l-1}^2+\vv^4h_{l-1})$ which is the same as in \cite{ahs15}. This means that the additional McKean-Vlasov component does not add computational complexity (per equation in the system of particles) and their conclusion about the computational cost of the method remains valid in our case. If $\delta \leq \e^2$ there is not benefit from using discretization methods customized for the small noise case. Moreover, if $\delta \geq e^{-\ff 1 \e}$, the EM scheme combined with the MLMC method leads to a cost $O(1).$ This is the same cost we would have with the standard MC method if we had $X_T$ as a formula of $W_T$, so no discretization method was required.

\end{document}